\RequirePackage{fix-cm}
\documentclass[smallextended]{svjour3}       % onecolumn (second format)
\smartqed  % flush right qed marks, e.g. at end of proof

\usepackage{graphicx}
\usepackage{amsmath}
\usepackage{amsfonts} 
\usepackage{natbib}

\newcommand{\bbR}{\mathbb{R}}

\newcommand{\Ic}{\mathcal{I}}
\newcommand{\argmax}{\mathop{\mathrm{argmax}}}

\newtheorem{assumption}{Assumption}
\bibpunct{(}{)}{;}{a}{}{,}
\begin{document}
\sloppy
\title{Submodular Optimization Problems and Greedy Strategies: A Survey%\thanks{Grants or other notes
%about the article that should go on the front page should be
%placed here. General acknowledgments should be placed at the end of the article.}
}
%\subtitle{Do you have a subtitle?\\ If so, write it here}

%\titlerunning{Short form of title}        % if too long for running head

\author{Yajing Liu \textsuperscript{1} \thanks{Yajing Liu\\yajing.liu@nrel.gov}\and
       Edwin K.~P. Chong\thanks{Edwin K.~P. Chong\\edwin.chong@colostate.edu}\textsuperscript{2} \and
Ali~Pezeshki\textsuperscript{2}\thanks{Ali~Pezeshki\\ali.pezeshki@colostate.edu}\and Zhenliang~Zhang\textsuperscript{3}\thanks{Zhenliang~Zhang\\zhenliang.zhang@alibaba-inc.com}
}

%\authorrunning{Short form of author list} % if too long for running head
\institute{\textsuperscript{1} National Renewable Energy Laboratory (NREL),
Golden, CO \\
         \textsuperscript{2} Department of Electrical and Computer Engineering, and Department of Mathematics,   Colorado State University, Fort Collins, CO\\
          \textsuperscript{3} Alibaba iDST, Seattle, WA
}

\date{Received: date / Accepted: date}
% The correct dates will be entered by the editor

\maketitle

\begin{abstract}
The greedy strategy is an approximation algorithm to solve optimization problems arising in decision making with multiple actions. How good is the greedy strategy compared to the optimal solution? In this survey, we mainly consider two classes of optimization problems where the objective function is submodular.  The first is set submodular optimization, which is to choose a set of actions to optimize a set submodular  objective function, and the second is string submodular optimization, which is to choose an ordered set of actions to optimize a string submodular function. Our emphasis here is on performance bounds for the greedy strategy in submodular optimization problems. Specifically, we review performance bounds for the greedy strategy, more general and improved bounds in terms of curvature,  performance bounds for the batched greedy strategy, and performance bounds for Nash equilibria.
% \PACS{PACS code1 \and PACS code2 \and more}
% \subclass{MSC code1 \and MSC code2 \and more}
\end{abstract}

\section{Introduction}
We are often faced with choosing a set of actions from a ground set of actions to optimize an objective function. Such problems arise in a multitude of applications of interest to discrete-event dynamic system researchers. A specific example is the task assignment problem  \citep{streeter2008, Zhang2016, Liu2018}, one of the fundamental combinatorial optimization problems in the study of optimization or operations research. This problem  involves a number of agents and a number of tasks. Each agent successfully accomplishes a task with a certain probability and the aim is to assign the available tasks to a given number of agents such that the probability of accomplishing the tasks is maximized.

When the number of agents is relatively small, we can use brute-force search  \citep{Bruteforce} to enumerate all possible candidate solutions to find the optimal solution. However, when the number of agents is large, it is impractical to enumerate all the possible candidate solutions. At this point, we have to resort to approximation methods.  One of the most well-studied approximation methods is  the \emph{greedy strategy} \citep{Nemhauser19781}, which starts  with the empty set and iteratively adds  to the current solution set an element that results in the largest gain in the objective function  while satisfying the constraints. The greedy strategy yields an approximation to an optimal solution in a reasonable amount of time. The downside is that there is often no theoretical guarantee for the greedy strategy. But when the problem has a special property called \emph{submodularity}, the greedy strategy is provably guaranteed to 
produce a solution with an objective value at least a constant scalar times the optimum value.  Celebrated results by  \cite{Fisher19782} and \cite{Nemhauser19781} prove that when the objective function $f$ is a monotone submodular set function with $f(\emptyset)=0$, the  greedy strategy yields a {$1/2$-approximation}\footnotemark
% ...
\footnotetext{The term $\beta$-approximation means that $f(G)/f(O)\geq \beta$, where $G$ and $O$ denote a greedy solution and an optimal solution, respectively.} for a general matroid and a $(1-e^{-1})$-approximation for a uniform matroid. 

For set optimization problems, the objective function is not influenced by the order of actions. However, a great number of problems in engineering and applied science aim to optimally choose a string (finite sequence) of actions over a finite horizon to maximize an objective function whose value depends on the order of actions. The problem arises in sequential decision making in engineering, economics, management science, and medicine. A motivating example is the problem of scheduling sensors to detect targets \citep{LiK09}. Suppose that a given number of sensors are distributed in a sensor field to detect a certain number of targets. The goal is to activate sensors sequentially to maximize the total coverage area. If the coverage region of each sensor remains constant over time, the total coverage area is not influenced by the order of the sensors activated, and the problem becomes a set optimization problem. However, if the sensors are moving, then the total coverage area depends on the order of the sensors activated, which makes the  problem fall into the framework of string optimization problems. The optimal solution to a string optimization problem is characterized by dynamic programming via Bellman's principle \citep{Powell2007}.  However, the approach  suffers from the curse of dimensionality and is therefore impractical for many problems of interest. This motivates the study of approximation algorithms, among which the greedy strategy is easy to implement and has guaranteed performance bounds under certain conditions. For example, \cite{streeter2008} prove that when the objective function is prefix and postfix monotone
and  has the diminishing-return property (as defined later in the paper), the greedy strategy yields a $(1-e^{-1})$-approximation.

In this paper, we review  the performance  guarantees for  greedy  strategies in submodular maximization problems. The paper is organized as follows. In Section~\ref{setactions}, we review results that are related to choosing sets of actions. This involves introducing set functions, set optimization problem, performance bounds for the greedy strategy, examples, curvature,  improved bounds, batched actions, and noncooperative games. In Section~\ref{stringaction}, we review results related to choosing strings of actions. This involves introducing new  notation and terminology, string optimization problem, performance of the greedy strategy, and applications. In Section~\ref{conclusions}, we conclude by listing a number of related papers that consider extensions and/or variation of greedy strategies and their performance bounds in combinatorial optimization problems.
\section{Sets of Actions}
\label{setactions}
In this section, we  first introduce our notation for sets, properties of set functions, and set optimization problems. Then, we review various performance bounds for the greedy strategy.
\subsection{Set Functions}
\label{setfunctions}
Before we introduce  functions defined  on sets, we would like to introduce some similar and familiar properties for functions defined on real numbers. Consider a real function $f:\mathbb{R}\rightarrow\mathbb{R}$. The function $f$ is said to be monotone and submodular if it satisfies properties 1 and 2 below, respectively:
\begin{itemize}
%\item WLOG, $f(0)=0$.
\item [1.] Monotone: $\forall x\leq y\in\mathbb{R}$, $f(x)\leq f(y)$.
\item [2.] Submodular: $\forall x\leq y\in\mathbb{R}$, $\forall z\in\mathbb{R}$, 
$f(x+z)-f(x)\geq f(y+z)-f(y)$
\end{itemize}
%\begin{figure}[ht]
%\begin{center}
%% Use the relevant command to insert your figure file.
%% For example, with the graphicx package use
%\vspace{-0.1in}
%\includegraphics[width=4in]{diminishing.pdf}
%% The figure caption is below the figure.
%\caption{Illustration of Submodularity}
%\end{center}
%\vspace{-0.25in}
%\end{figure}

\begin{figure}
   \centering
   \includegraphics[width=\columnwidth, trim= 50 120 100 60, clip]{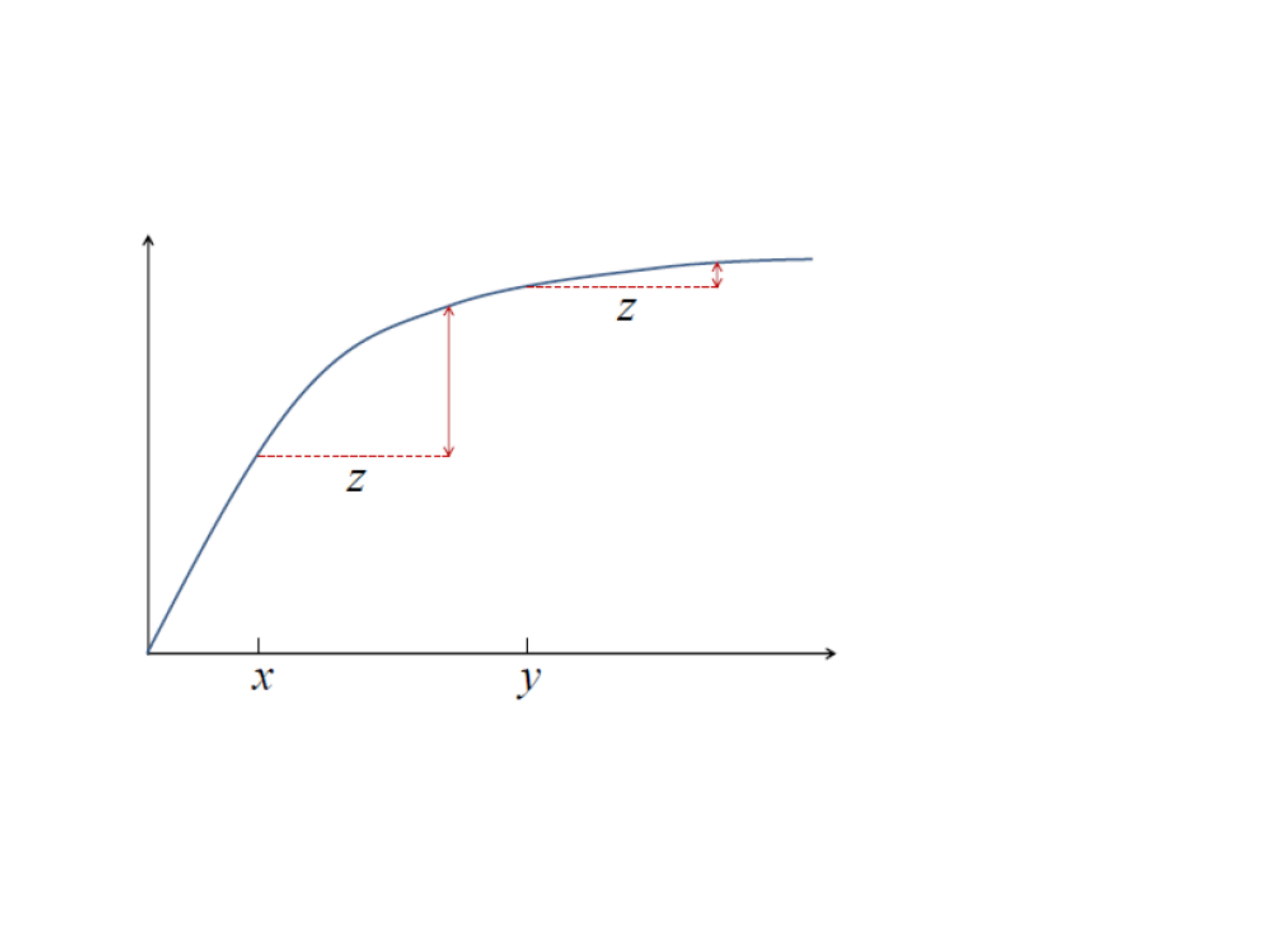}
   \caption{Illustration of submodularity}
\label{Fig2}
\end{figure}

The `monotone' property here simply means being `nondecreasing'.  The function in Fig.~1 satisfies the monotone property. 
From Fig.~1, we can see that the function is a concave function -- adding $z$ to $x$ gains more than adding $z$ to $y$, which tells us that the additional value accrued by  adding a  number to a smaller number is larger than adding it to a bigger number. This is consistent with the inequality $f(x+z)-f(x)\geq f(y+z)-f(y)$ for $x\leq y$, so we say that `submodularity' here boils down to `concavity' in some sense. 

In this paper, we want to go beyond the real line to a more general setting. Specifically, we will consider objective functions with multiple decision ``actions" as arguments. The first setting is sets of actions, and the second one is strings (ordered sets) of actions. We introduce functions defined on sets first.

Let $X$ denote a ground set, which includes all possible actions. Let $2^X$ denote the power set of $X$, which includes all possible subsets of $X$. The \emph{size} or \emph{cardinality} of a set $S\in 2^X$ is denoted by $|S|$, and the \emph{empty set} is denoted by $\emptyset$. Define a set function $f$: $2^X\longrightarrow \mathbb{R}$.
The set function $f$ is said to be \emph{monotone} and \emph{submodular} if  it satisfies properties~i and~ii below, respectively:
\begin{itemize}
\item [i.] Monotone: $\forall A\subseteq B\subseteq X$, $f(A)\leq f(B)$.
\item [ii.] Submodular: $\forall A\subseteq B\subseteq X$ and $\forall j\in X\setminus B$, $f(A\cup\{j\})-f(A)\geq f(B\cup\{j\})-f(B)$.
\end{itemize}
Notice the similarity between these properties and those involving functions on the real line introduced earlier.

For convenience, we denote the incremental value of adding a set $T$ to the set $A\subseteq X$ as $\varrho_T(A)=f(A\cup T)-f(A)$ (following the notation in \cite{conforti1984submodular}).

A set function $f$: $2^X\longrightarrow\mathbb{R}$ is called a \emph{polymatroid set function} \citep{Boros2003} if it is monotone, submodular, and $f(\emptyset)=0$.
Submodularity in property~ii means that the additional value accruing from an extra action decreases as the size of the input set increases, and is also called the \emph{diminishing-return} property in economics.
%Submodularity implies that for any $A\subseteq B\subseteq X$ and $T\subseteq X\setminus B$, 
%\begin{equation}
%\label{eqn:submodularimplies}
%f(A\cup T)-f(A)\geq f(B\cup T)-f(B).
%\end{equation}
Submodularity has many equivalent definitions; for example, $f:2^X\longrightarrow\mathbb{R}$ is submodular if $\forall A, B\subseteq X$, $f(A)+f(B)\geq f(A\cup B)+f(A\cap B)$. For more equivalent definitions, see \cite{Nemhauser19781}. 

The set function $f$ is called \emph{supermodular} if $-f$ is submodular. Moreover, $f$ is called \emph{modular} if it is both submodular and supermodular, i.e., for any $A\subseteq B\subseteq X$,
\begin{equation}
\label{modular}
f(A)+f(B)= f(A\cup B)+f(A\cap B).
\end{equation}
By induction, (\ref{modular}) implies that for any $S\subseteq X$,   
\begin{equation}
\label{modular1}
f(S)-f(\emptyset)=\sum\limits_{s\in S}(f(\{s\})-f(\emptyset)).
\end{equation}
 By (\ref{modular1}),  $f-f(\emptyset)$ is additive when $f$ is modular.
If $f(\emptyset)=0$, then $f(S)=\sum_{s\in S}f(\{s\})$, which implies that  $f$ is additive. It is also easy to check that $f$ is modular iff for any subset $S\subseteq X$, 
\begin{equation}
\label{modular2}
f(S)=\omega(\emptyset)+\sum_{s\in S}\omega(s)
\end{equation}
for some weight function $\omega: X \rightarrow \mathbb{R}$ \citep{Krause2013}.

There are many non-trivial examples of submodular or supermodular set functions. We only consider submodular maximization problems in this paper, so we only give submodular function examples. For supermodular examples, see \cite{Lovasz1983}. To easily understand submodularity, we provide a simple example as follows.
\begin{example}
\textbf{Sensor Coverage.} Let $X$ be a family of locations in space where we can place sensors. If a sensor  is placed at a particular location in space, it covers a circular area around it as illustrated in Fig.~\ref{Fig2}. Let $f(S)$ denote the total area covered if we place sensors at locations $S\subseteq X$ (see Fig.~\ref{Fig2}). The set function $f$ is submodular. An instance of submodularity is illustrated in the figure. As can be seen, the gain in adding sensor $3$ after placing sensor $1$ is larger than the gain in adding sensor $3$ after placing sensors $1,2$. $\qed$
\end{example}
\begin{figure}
   \centering
   \includegraphics[width=\columnwidth, trim= 100 180 100 20, clip]{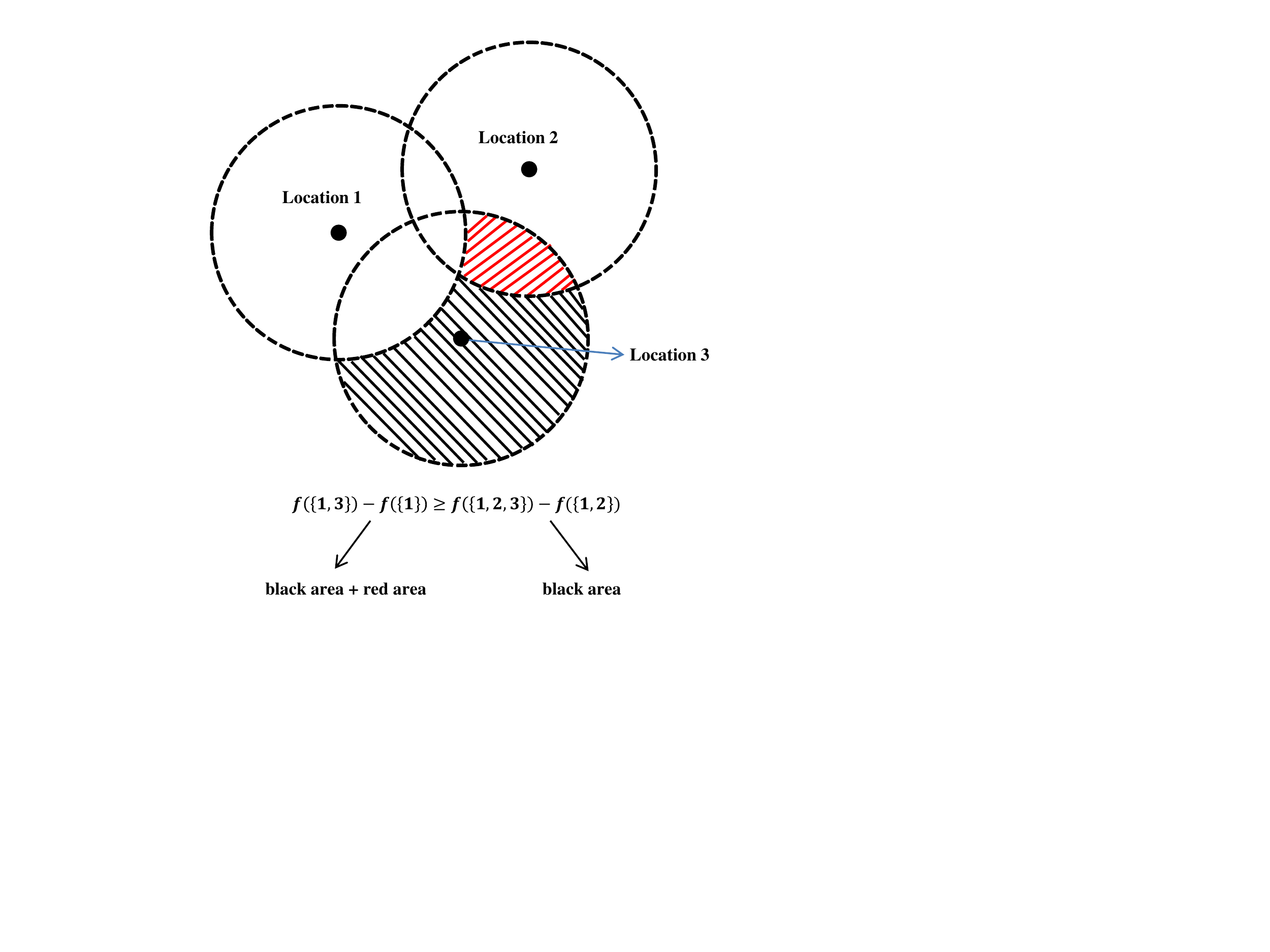}
   \caption{Sensor coverage as an example of a submodular function}
\label{Fig2}
\end{figure}
Submodular functions arise in many applications,  such as the rank  function of the matrix formed by its columns, weighted coverage functions, the rank function of a matroid, Shannon entropy, mutual information, cut capacity functions, some measurements on the graph, etc. \citep{Lovasz1983,Krause2013}
\subsection{Submodular Set Optimization Problem}
\label{setoptimization}
Submodular set optimization plays an important role in combinatorial optimization. It has a wide range of applications, including  generalized assignment \citep{Shmoys1993,Cohen2006,Nauss2003,Fleischer2006,Bator1957,Korula2015,Vondrak2008}, matroid partition \citep{Edmonds1965,Cunningham1986,Knuth1973}, maximum cut \citep{Goemans1995,Sahni1976}, maximum coverage location  \citep{Church1974,Khuller1999,Fisher1977}, multi-agent coverage problem \citep{Sun2017}, leader-selection problem in multi-agent systems \citep{Clark2011}, welfare maximization \citep{Korula2015,Vondrak2008,Kapralov2013}, and data summarization \citep{Lin2011, Badanidiyuru2014, mirzasoleiman17a}.  The aim is to find a set of actions satisfying some constraints to maximize the objective function. The set optimization problem can be formulated as follows:
\begin{align}\label{setproblem}
\begin{array}{l}
\text{maximize} \ \    f(M), \ \quad \text{subject to} \ \  M\in \mathcal{I},
\end{array}
\end{align}
where $\mathcal{I}$ is a non-empty collection of subsets of a finite set $X$, and $f$ is a real-valued  submodular set function defined on the power set $2^X$ of $X$. Before proceeding any further with discussing optimization problem (\ref{setproblem}), we will need to introduce some concepts related to the constraint set $\mathcal{I}$.

Let $X$ be a finite set, and $\mathcal{I}$ be a non-empty collection of subsets of  $X$. The collection $\mathcal{I}$ is said to be \emph{hereditary} if it satisfies property i below and  has the \emph{augmentation} property if it satisfies property ii below:
\begin{itemize}
\item  [i.] Hereditary: For all $B\in\mathcal{I}$, any set $A\subseteq B$ is also in $\mathcal{I}$.
\item  [ii.] Augmentation: For any $A,B\in \mathcal{I}$, if $|B|>|A|$, then there exists $j\in B\setminus A$ such that $A\cup\{j\}\in\mathcal{I}$.
\end{itemize}
The pair $(X,\mathcal{I})$ is called an \emph{independence system} if it satisfies property~i. In this case, the sets in $\mathcal{I}$ are called \emph{independent sets}. A \emph{maximal independent set} is an independent set that is not a subset of any other independent set \citep{conforti1984submodular}. The independence system  $(X,\mathcal{I})$ is called a \emph{matroid}   if it satisfies  property~ii \citep{Edmonds2003}. The pair $(X,\mathcal{I})$ is called a \emph{uniform matroid} if $\mathcal{I}=\{S\subseteq X: |S|\leq K\}$ for a given $K$ \citep{Nemhauser19781}.  All maximal independent sets in a  matroid have the same cardinality. We call this cardinality the \emph{rank} of the matroid. In the uniform matroid above, the rank is $K$.
\begin{example}
\label{independenceexample}
 We now give three example collections to illustrate the notions of independence systems and matriods.
 Let $X = \{a, b, c\}$, $\mathcal{I}_1= \{\{a\}, \{b\}, \{a,c\}, \{c\}, \emptyset\}$, 
$\mathcal{I}_2= \{\{a\}, \{a,b\}\}$, and $\mathcal{I}_3=\{\emptyset, \{a\},\{b\},\{a,b\}\}$. It is easy to check that $\mathcal{I}_1$ satisfies the hereditary property but not augmentation, $\mathcal{I}_2$ satisfies augmentation but not the hereditary property, and $\mathcal{I}_3$ satisfies both hereditary  and augmentation properties. Hence, $(X,\mathcal{I}_1)$ is an independence system, $(X,\mathcal{I}_3)$ is a matroid, and $(X,\mathcal{I}_2)$ is neither an independence system nor a matroid. The maximal independent sets in $(X,\mathcal{I}_1)$ are $\{b\}$ and $\{a, c\}$, and $(X,\mathcal{I}_3)$ only has one maximal independent set $\{a, b\}$. $\qed$
\end{example}

Let $(X,\mathcal{I})$ be an independence system where $\mathcal{I}$ is nonempty, and let $S\subseteq X$ be an arbitrary subset of $X$.  A \emph{basis} of $S$ is  a  subset $B$ of $S$  that satisfies the following two conditions:
1. It is an independent set; i.e., $B\in\mathcal{I}$.
2. It is maximal; i.e., $B$ is not a subset of any other independent subset of $S$. The subset $B$ satisfying the above two conditions is also called a maximal independent subset of $S$.
Define
\[\mbox{\emph{lower  rank}  of} \ S=\text{lr}(S)=\min\{|B|: B\ \mbox{is a   basis of}\ S\},\]
\[\mbox{\emph{upper  rank}  of} \ S=\text{ur}(S)=\max\{|B|: B\ \mbox{is  a basis  of}\ S\}.\]
Note that lr($S$) and ur($S$) might not be well defined, depending on $S$. Note also that in the definition above, $S$ is not necessarily in $\mathcal{I}$.
The number 
\begin{equation}
\label{rankquotientdef}
q(X,\mathcal{I})=\min\left\{\frac{\text{lr}(S)}{\text{ur}(S)}:S\subseteq X \ \text{and} \ \text{ur}(S)> 0\right\}
\end{equation}
is called the \emph{rank quotient} of $(X,\mathcal{I})$   \citep{hausmann1980}.

\begin{example}
 To illustrate the concept of rank quotient, again consider the independence system $(X,\Ic_1)$ given in Example~\ref{independenceexample}. We now consider all the subsets of $X$ and calculate their lower and upper ranks. If $S$ is a singleton (i.e., $\{a\}$, $\{b\}$, or $\{c\}$), then $S$ has only one basis, which is $S$ itself. In this case, $\text{lr}(S)=\text{ur}(S)=1$, which means that $\text{lr}(S)/\text{ur}(S)=1$. 
 
 If $S=\{a,b\}$, its bases are $\{a\}$ and $\{b\}$. Again, $\text{lr}(S)=\text{ur}(S)=1$, which means that $\text{lr}(S)/\text{ur}(S)=1$. Note that $\{a,b\}$ is not a basis of $S$ because it does not belong to $\Ic_1$. 
If $S=\{a,c\}$, it has only one basis, which is itself, and again $\text{lr}(S)/\text{ur}(S)=1$. If $S=\{b,c\}$, its bases are $\{b\}$ and $\{c\}$, in which case $\text{lr}(S)/\text{ur}(S)=1$ again.

If $S=\{a,b,c\}=X$, the bases are $\{b\}$ and $\{a,c\}$. So, $\text{lr}(S)=1$ and $\text{ur}(S)=2$, which implies that $\text{lr}(S)/\text{ur}(S)=1/2$.

Because the rank quotient is the smallest among the ratios calculated above, we deduce that $q(X,\Ic_1)=1/2$. $\qed$
\label{rankquotient}
\end{example}

\begin{example}
 As in Example~\ref{rankquotient}, we can similarly check that $q(X,\Ic_3)=1$. 
In fact, the rank quotient of any matroid $(X,\Ic)$ is equal to $1$, because for any susbset $S\subseteq X$, $\text{lr}(S)=\text{ur}(S)$ \citep{Edmonds1966}. The rank quotient of an independence system  $(X,\mathcal{I})$ can be regarded as a measure of how much $(X,\mathcal{I})$ differs from being a matroid.
 $\qed$

\end{example}

For any independence system $(X,\mathcal{I})$, if there exist matroids 
 $(X,\mathcal{I}^i)$ ($1\leq i\leq p$) such that $\mathcal{I}=\mathcal{I}^1\cap\cdots \cap \mathcal{I}^p$, then the pair $(X,\mathcal{I})$ is called the \emph{intersection} of the matroids $(X,\mathcal{I}^i)$ \citep{hausmann1980}.

Finding the optimal solution to  (\ref{setproblem}) in general is NP-hard. The \emph{greedy strategy} provides a tractable way to approximately solve the problem, which starts with the empty set, and incrementally adds an element to the current solution set giving the largest gain in the objective function under the constraints. Although the greedy strategy yields an approximate solution, its performance might be arbitrarily poor. However, when the optimization problem has the further special structure of being polymatroid, the greedy strategy has provable guarantees. The celebrated results by  \protect\cite{Fisher19782} and \protect\cite{Nemhauser19781} show that the greedy strategy provides a good approximation to the optimal solution when the objective function is a polymatroid set function under both general matroid constraints and uniform matroid constraints. We will review the performance of the greedy strategy for (\ref{setproblem}) under different constraints in the following section. 
\subsection{Performance Bounds for Greedy Strategy}
\label{greedybounds}
First we introduce definitions of the optimal strategy and the greedy strategy.
\noindent\textbf{Optimal Set}:
Any set $O$ is called an \emph{optimal solution} of Problem~(\ref{setproblem}) if
\[O\in \argmax_{{M\in\mathcal{I}}}f(M),\]
where argmax denotes the set of actions that maximize $f(\cdot)$.

\noindent\textbf{Greedy Algorithm}: 

\noindent\textbf{Input}: A pair $(X,\mathcal{I})$, a set function $f:2^X\rightarrow \mathbb{R}$

\noindent\textbf{Output}: A subset $G\in\mathcal{I}$

\noindent$G_0\leftarrow \emptyset$

\noindent For $i=1,2,\ldots$,

\noindent  $g_i\leftarrow \argmax\limits_{\substack{a\in X\setminus G_{i-1} \\ G_{i-1}\cup \{a\}\in\mathcal{I}}}f(G_{i-1}\cup a)$. If $g_i\neq \emptyset$, set $G_i = G_{i-1}\cup\{g_i\}$; otherwise, stop and set $G = G_{i-1}$.

Any output of the above algorithm is called a \emph{greedy solution}.
Note that there may exist more than one optimal solution or more than one greedy solution.
How good is a  greedy solution compared to an optimal solution in terms of the objective function? In the following theorems, we review performance bounds for the greedy strategy under different constraints. These bounds are worst-case performance bounds, which means that the greedy strategy performs much better than those bounds in many cases.

\begin{theorem}\citep{hausmann1980}
\label{independencethm}
Let $(X,\mathcal{I})$ be an independence system. If $f$ is additive on $X$, i.e., $f(S)=\sum_{s\in S} f(\{s\})$ for any subset $S\subseteq X$, then any greedy solution $G$ satisfies
\begin{equation}
\label{independencebound}
\frac{f(G)}{f(O)}\geq q(X,\mathcal{I}),
\end{equation}
\end{theorem}
where $q(X,\mathcal{I})$ is the rank quotient defined in Section~\ref{setoptimization}. Furthermore, for some function $f$, (\ref{independencebound}) holds with equality.

\begin{remark}
\label{remark2}
When $(X,\mathcal{I})$ is a matroid, $q(X,\mathcal{I})=1$. By Theorem~\ref{independencethm}, the  greedy strategy is optimal when $(X,\mathcal{I})$ is a matroid and the objective function is additive.
\end{remark}
\begin{remark}
\label{remark3}
When $(X,\mathcal{I})$ is the intersection of $p$ matroids, then $q(X,\mathcal{I})\geq 1/p$ \citep{hausmann1980}. So when $p=1$, i.e., $(X,\mathcal{I})$ is a matroid, the  greedy strategy  is optimal, which is consistent with Remark~\ref{remark2}.
\end{remark}
\begin{example}
We provide an example\footnotemark
% ...
\footnotetext{We thank the anonymous reviewer for this example.} to demonstrate the performance bound in Theorem~\ref{independencethm}. Let $X  = \{s,t, u, v, w, x\}$, and consider the collection of subsets  \[\begin{aligned}
    \mathcal{I} =\{&\emptyset, \{s\}, \{t\}, \{s,t\}, \{u\}, \{v\}, \{w\}, \{x\}, \{u,v\}, \{u,w\},  \{u,x\}, \{v,w\}, \{v,x\}, \\ &\{w,x\},
    \{u,v,w\}, \{u,v,x\}, \{u,w, x\}, \{v,w, x\}, \{u,v, w, x\} \}.   
\end{aligned}\]
Define a function $f$ such that $f(A)=\sum_{a\in A}f(\{a\})$. Let $f(\{s\})=1.01, f(\{u\})=f(\{v\})=f(\{w\})=f(\{x\})=1$, and $f(\{t\})=0$.

It is easy to check that $(X,\mathcal{I})$ is an independence system. If $S=\{s,u,v,w,x\}$, it has bases $\{s\}$ and $\{u,v,w,x\}$, which results in $\text{lr}(S)/\text{ur}(S)=1/4$.
Because the maximum cardinality of the maximal independent subsets of any $S\subseteq X$  is 4, $\text{lr}(S)/\text{\text{ur}(S)}\geq {1}/{4}$ for any set $S\subseteq X$ with $\text{ur}(S)>0$. Therefore, $q(X,\mathcal{I})={1}/{4}$. The greedy solution is $G=\{s,t\}$ with $f(G)=1.01$ and the optimal solution is $O=\{u,v,w,x\}$ with $f(O)=4$, which satisfy the bound $f(G)/f(O)\geq q(X,\mathcal{I})$. In fact, the bound holds with equality if we lower $f(\{s\})$ to exactly $1$. $\qed$
 
\end{example}

The following theorem bounds the performance of the greedy strategy  when $(X,\mathcal{I})$ is the intersection of $p$ matroids and $f$ is a polymatroid set function.

\begin{theorem}\citep{Fisher19782}
\label{intersectionmatroidthm}
Let $(X,\mathcal{I})$ be the intersection of $p$ matroids and  $f:2^X\rightarrow\mathbb{R}$  a polymatroid set function. Then any greedy solution $G$ satisfies
\begin{equation}
\label{intersectionmatroidbound}
\frac{f(G)}{f(O)}\geq \frac{1}{1+p}.
\end{equation}
\end{theorem}
\begin{remark}
\label{remark4}
The condition that $f$ is additive in Theorem~\ref{independencethm} is stronger than the condition that $f$ is a polymatroid set function in Theorem~\ref{intersectionmatroidthm}, so the bound $1/p$ in Theorem~\ref{independencethm}
is stronger than the bound $1/(1+p)$ in Theorem~\ref{intersectionmatroidthm}.

\end{remark}
\begin{remark}
The bound $1/(1+p)$ can be achieved for any positive integer $p$. When $p=1$, $(X,\mathcal{I})$ is a matroid, and the bound becomes $1/2$, which means that the greedy strategy yields $1/2$-approximation for general matroid constraints.
\end{remark}

\begin{remark}
Theorem \ref{intersectionmatroidthm} requires that $f(\emptyset)=0$. If $f(\emptyset)\neq 0$, the following performance bound holds
 \[\frac{f(G)-f(\emptyset)}{f(O)-f(\emptyset)}\geq \frac{1}{1+p}.\]

\end{remark}

The following theorem provides a performance bound for the greedy strategy when $(X,\mathcal{I})$ is a uniform matroid and $f$ is a polymatroid set function.
\begin{theorem}\citep{Nemhauser19781}
\label{unifrommatroidboundthm}
Let $(X,\mathcal{I})$ be a uniform matroid and  $f:2^X\rightarrow\mathbb{R}$  a polymatroid set function. Then any greedy solution $G_K$ satisfies 
\begin{equation}
\label{intersectionmatroidbound}
\frac{f(G)}{f(O)}\geq 1-\left(1-\frac{1}{K}\right)^K>1-\frac{1}{e},
\end{equation}
where $K$ is the rank of the uniform matroid and $e$ is the base of the natural logarithm.
\end{theorem}
\begin{remark}
The bound $1-(1-1/K)^K$ is stronger than the bound $1/(1+p)$ when $p=1$ in Theorem~\ref{intersectionmatroidthm}, because the uniform matroid is a special matroid. 
\end{remark}
\begin{remark}
The bound $1-(1-1/K)^K$ is decreasing in $K$ and tends to $1-1/e$ when $K$ goes to infinity. When $K=1$, the bound becomes 1, which is consistent with the fact that the greedy strategy chooses the best action at each stage.
\end{remark}
\begin{remark}
The bound $1-(1-1/K)^K$  is tight, which means that it can be achieved for each $K$ \citep{Nemhauser19781}.
\end{remark}
\begin{remark}
By Theorem~\ref{intersectionmatroidthm}, the greedy strategy only achieves a $1/2$-approximation under general matroid constraints. However, \cite{Calinescu2011} proved that a variant of the greedy strategy yields a $(1-1/e)$-approximation under general matroid constraints.

\end{remark}
\subsection{Examples}
\label{examples}
We introduce two examples -- a task scheduling problem and an adaptive sensing problem -- to illustrate  polymatroid set functions. In both problems, $(X,\mathcal{I})$ is a uniform matroid and hence the greedy strategy satisfies a $(1-e^{-1})$-approximation.

\textbf{Task Assignment Problem:}
The task scheduling problem  was posed by \cite{streeter2008}, and was also analyzed in \cite{Zhang2016} and \cite{Liu2018}.  In this problem, there are $n$ subtasks and a  set $X$ of $N$ agents. At each stage, a subtask $i$ is assigned to an agent $a$, who accomplishes the task with probability $p_i(a)$. Let $X_i(\{a_1,a_2,\ldots, a_k\})$ denote the Bernoulli random variable that signifies whether or not subtask $i$ has been accomplished after assigning the set of agents $\{a_1,a_2,\ldots, a_k\}$ over $k$ stages. Then $\frac{1}{n}\sum_{i=1}^n	X_i(\{a_1,a_2,\ldots,a_k\})$ is the fraction of subtasks accomplished after $k$ stages by employing agents $\{a_1,a_2,\ldots, a_k\}$. The objective function $f$ for this problem is the expected value of this fraction, which can be written as
\begin{equation}
\label{objectivefunction}
f(\{a_1,\ldots,a_k\})=\frac{1}{n}\sum_{i=1}^n\left(1-\prod_{j=1}^k\left(1-p_i(a_j)\right)\right).
\end{equation}
The aim is to choose a set of agents to maximize this objective function.

Assume that $p_i(a)>0$ for any $a\in X$. Then it is easy to check that $f$ is monotone, submodular, and $f(\emptyset)=0$, which implies that $f$ is a polymatroid set function. 

\textbf{Adaptive Sensing:}
As our second example application, we consider the adaptive sensing design problem posed in \cite{Zhang2016} and \cite{Liu2018}. Consider a signal of interest $x \in {\rm I\!R}^2$ with normal prior distribution $\mathcal{N} (0, I)$, where $I$ is the $2\times 2$ identity matrix; our analysis easily generalizes to dimensions larger than $2$. Let $\mathbb{B}=\{\mathrm{Diag}(\sqrt{b},\sqrt{1-b}):  b\in\{b_1,\ldots,b_N\}\}$, where $ b_i\in[0,1]$ for $1\leq i
\leq N$. At each stage $i$, we make a measurement $y_i$ of the form
\[
y_i=B_ix+w_i, 
\]
where $B_i \in\mathbb{B}$ and $w_i$ represents i.i.d.\ Gaussian measurement noise with mean zero and covariance $\sigma^2 I$, independent of $x$.

The objective function $f$ for this problem is the information gain, which can be written as
\[
f(\{B_1,\ldots, B_k\})=H_0-H_k.
\]
Here, $H_0=\frac{N}{2}\text{log}(2\pi e)$ is the entropy of the prior
distribution of $x$ and $H_k$ is the entropy of the posterior
distribution of $x$ given $\{y_i\}_{i=1}^k$; that is,
\[
H_k=\frac{1}{2}\text{log det}(P_k)+\frac{N}{2}\text{log}(2\pi e),
\] 
where for $1\leq j\leq k$
\[
P_j=\left(P_{j-1}^{-1}+\frac{1}{\sigma^2}B_j^TB_j\right)^{-1}
\]
is the posterior covariance of $x$ given $\{y_i\}_{i=1}^j$. The objective is to choose a set of measurements to maximize the information gain $
f(\{B_1,\ldots, B_K\})=H_0-H_K$. 

It is easy to check that $f$ is monotone, submodular, and $f(\emptyset)=0$; i.e., $f$ is a polymatroid set function.
\subsection{Curvature}
\label{curvature}
As we saw in Section~\ref{setfunctions}, submodularity is a \emph{second-order} property by analogy to concavity. If we can \emph{quantify} this second order property, then we can get tighter bounds. One way to quantify the second order property is to use the \emph{total curvature}, defined by \cite{conforti1984submodular}:
\[c(f):=\max_{j\in X: \varrho_j(\emptyset)\neq 0}\left\{1-\frac{\varrho_j({X\setminus\{j\}})}{\varrho_j(\emptyset)}\right\}.\]
To see that this is a second-order property, rewrite it in terms of differences of differences:
\[c(f):=\max_{j\in X: f(\{j\})\neq f(\emptyset)}\left\{\frac{(f(\{j\})-f(\emptyset))-(f(X)-f(X\setminus \{j\}))}{f(\{j\})-f(\emptyset)}\right\}.\]
For convenience, we use $c$ to denote $c(f)$ when there is no ambiguity.
Note that $0\leq c\leq 1$ when $f$ is a polymatroid set function, and $c=0$ when $f$ is {modular}. When $f$ is modular, $f-f(\emptyset)$ is additive. If we consider $f-f(\emptyset)$ as the objective function, then the greedy strategy achieves optimality.  Therefore, in the rest of the paper, when we assume that $f$ is a polymatroid set function, we only consider $c\in(0,1]$.

\cite{conforti1984submodular} provided performance bounds in terms of the total curvature for the greedy strategy under independence system, general matroid, and uniform matroid constraints, which will be reviewed as follows.
\begin{theorem}\citep{conforti1984submodular}
If $(X,\mathcal{I})$ is an independence system with ur$(X)=K$ and lr$(X)=k$, and $f$ is  is a polymatroid set function with total curvature $c$, then any greedy solution $G_K$ satisfies 
\[\frac{f(G_K)}{f(O)}\geq\frac{1}{c}\left[1-\left(1-\frac{c}{K}\right)^k\right],\]
and this bound is tight for all $0< c\leq 1$.
\end{theorem}

\begin{theorem}\citep{conforti1984submodular}
\label{greedyboundcurvature}
If $(X,\mathcal{I})$ is a matroid and $f$ is a polymatroid set function with total curvature $c$, then any greedy solution $G_K$ satisfies 
\[\frac{f(G_K)}{f(O)}\geq \frac{1}{1+c}.\]
Moreover, if $(X,\mathcal{I})$ is a uniform matroid with rank $K$, then 
any greedy solution $G_K$ satisfies 
\[\frac{f(G_K)}{f(O)}\geq \frac{1}{c}\left[1-\left(1-\frac{c}{K}\right)^K\right]>\frac{1-e^{-c}}{c}.\]
\end{theorem}

\begin{remark}
When $(X,\mathcal{I})$ is a matroid, the bound $1/(1+c)$ is stronger than the bound $1/2$ in Theorem~\ref{intersectionmatroidthm} because $c\in (0,1]$ when $f$ is a polymatroid set function and $1/(1+c)$ is nonincreasing in $c$. 
\end{remark}
\begin{remark}
The function $(1-e^{-c})/c$ is nonincreasing in $c$, and therefore $(1-e^{-c})/c\in [1-e^{-1}, 1)$ when $f$ is a polymatroid set function. Also it is easy to check that $(1-e^{-c})/c\geq 1/(1+c)$ for $c\in(0,1]$, which implies that the bound $(1-e^{-c})/c$ for the uniform matroid constraints is stronger than the bound $1/(1+c)$ for the general matroid constraints.
\end{remark}
\begin{remark}
The two bounds in terms of the total curvature $c$ are both tight; for proofs, see \cite{conforti1984submodular}.
\end{remark}
\begin{remark}
There are other notions of curvatures that can be used to characterize the second-order property of the set function $f$, such as the greedy curvature defined by \cite{conforti1984submodular} and the elemental curvature defined by \cite{Wang2014}. Performance bounds in terms of the corresponding curvatures  were also derived by \cite{conforti1984submodular} and  \cite{Wang2014} under different constraints.
\end{remark}

\textbf{Example:} Consider again the task assignment example from Section~\ref{examples}. For convenience, we only consider the special case $n=1$; our analysis can be generalized to any $n\geq 2$. For $n=1$, we have 
$$f(\{a_1,\ldots,a_k\})=1-\prod_{j=1}^k\left(1-p(a_j)\right),$$
where $p(\cdot)=p_1(\cdot)$.

Let us  order the elements of $X$ as $a_{[1]}, a_{[2]}, \ldots, a_{[N]}$ such that $$0<p(a_{[1]})\leq p(a_{[2]})\leq\ldots\leq p(a_{[N]})\leq 1.$$ Then by the definition of the total curvature $c$, we have 
\[
c=\max\limits_{j\in {X}}\left\{1-\frac{f(X)-f(X\setminus\{j\})}{f(\{j\})-f(\emptyset)}\right\}=1-\prod_{l=2}^N(1-p(a_{[l]})) < 1,
\]
which is consistent with our conclusion that $c\in [0,1]$.

\subsection{Improved Bounds}
\label{improvedbounds}
The performance bounds of \cite{conforti1984submodular} reviewed in Section~\ref{curvature}, are the best bounds in terms of the total curvature $c$  for general matroid constraints and uniform matroid constraints, respectively.  However, the total curvature $c$ depends on function values on sets outside the constraint matroid. If we are given a function defined only on the matroid,  problem (\ref{setproblem}) still makes sense, but the bounds involving $c$ do not apply. \cite{Liu20181,Liu20173} investigated modified bounds that overcome this drawback. The idea is first to extend a polymatroid set function defined on the matroid to one defined on the entire power set, and then apply the results from  \cite{conforti1984submodular}. However, not every polymatroid function defined on the matroid can be extended to one defined on the entire power set.

In \cite{Liu20173}, they first provide necessary and sufficient conditions for the existence of an incremental extension of a polymatroid set function defined on the uniform matroid of rank $k$ to one defined on the uniform matroid of rank $k+1$. Whenever a polymatroid objective function defined on a matroid can be extended to the entire power set, the greedy approximation bounds involving the total curvature of the extension apply. However, the bounds still depend on sets outside the matroid.  Motivated by this, \cite{Liu20173}  defined a new notion of curvature called \emph{partial curvature}, involving only sets in the matroid. They derived necessary and sufficient conditions for an extension of the function to have a total curvature that is equal to the partial curvature. Moreover, they proved that the bounds in terms of the partial curvature are in general improved over the previous ones.

The following theorems state the necessary and sufficient conditions for the existence of an extension of a polymatroid set function defined on the uniform matroid of rank $k$ to one defined on the uniform matroid of rank $k+1$.

\begin{theorem}\citep{Liu20173}
\label{polymatroidextension}
Let  $f:\mathcal{I}\rightarrow \mathbb{R}$ be a polymatroid  function defined on the uniform matroid  of rank $k$.  Then $f$ can be extended to a  polymatroid function $g$ defined on the uniform matroid of rank $k+1$  if and only if  for any  $A\subseteq X$ with $|A|=k+1$, any  $B\subset A$ with $|B|=k
$, and any $a\in B$, 
\begin{equation}
\label{sufneccond1}
f(B)-f(B\setminus\{a\})\geq f(B^*)-f(A\setminus\{a\}),
\end{equation}
where $B^*\in\argmax\limits_{\substack{B: B\subset A, |B|=k}}f(B)$.
\end{theorem}

\noindent\textbf{Construction}:  If $f$ is extendable, then an extension $g$ can be constructed as follows: 
For any $A$ with $|A|\leq k$, $g(A)=f(A)$;
For any $A$ with $|A|=k+1$, 
\begin{equation}
\label{construction}
g(A)=
g(B^*)+d_A,
\end{equation}
where $d_A$ satisfies 
\begin{equation}
\label{dsatisfy}
0\le d_A\le \min\limits_{\substack{B:B\subset A,  |B|=k \\a:a\in B}}\{f(B)-f(B^*)+f(A\setminus\{a\})-f(B\setminus \{a\})\}.
\end{equation}
\newcommand{\Ec}{\mathcal{E}}
Note that  $f:2^X\rightarrow \mathbb{R}$ is itself an extension of $f$ from $\mathcal{I}$ to the entire $2^X$, and the extended $f:2^X\rightarrow \mathbb{R}$ is a polymatroid function on $2^X$. Therefore, we have that $c(f)\geq d=\inf_{g\in \Ec_f} c(g)$, where $\Ec_f$ is the set of all polymatroid functions $g$ on $2^X$ that agree with $f$ on $\mathcal{I}$. So  if a polymatroid set function defined on the matroid can be extended to one defined on the whole power set, applying the performance bounds in Theorem~\ref{greedyboundcurvature} results in the following theorem.

\begin{theorem}\citep{Liu20173}
\label{Improvedbound}
 Let $(X,\mathcal{I})$ be a matroid of rank $K$ and $f:\mathcal{I}\rightarrow\mathbb{R}$ a polymatroid function.  If there exists an extension of $f$ to the entire power set, then any greedy solution $G_K$ to problem~$(\ref{setproblem})$ satisfies
 \begin{equation}
 \label{bound1}
 \frac{f(G_K)}{f(O)}\geq \frac{1}{1+d},
 \end{equation}
where $d=\inf_{g\in \Ec_f} c(g)$. In particular, when  $(X,\mathcal{I})$ is a uniform matroid, any greedy solution $G_K$ to problem~$(\ref{setproblem})$ satisfies
\begin{align}
\label{bound2}
\frac{f(G_K)}{f(O)}&\geq\frac{1}{d}\left[1-\left(1-\frac{d}{K}\right)^K\right]> \frac{1}{d}\left(1-e^{-d}\right).
\end{align}
\end{theorem}
\begin{remark}
The bounds $1/(1+d)$ and $(1-e^{-d})/d$  apply to  problems where the objective function is a polymatroid function defined only for sets in the matroid and can be extended to one defined on the entire power set. However, these bounds still depend on sets not in the matroid, because of the way $d$ is defined.
\end{remark}

Then \cite{Liu20173} defined a new curvature called the \emph{partial curvature} $b(h)$ as follows:
\begin{equation}
\label{curvature_b}
b(h):=\max_{\substack{j, A: j\in A\in\mathcal{I}\\ h(\{j\})\neq h(\emptyset)}}\left\{1-\frac{h(A)-h(A\setminus\{j\})}{h(\{j\})-h(\emptyset)}\right\},
\end{equation}
and the partial curvature satisfies that $b(f)\leq c(g)$ when $g$ is an extension of $f$ from $\mathcal{I}$ to $2^X$. The following theorem provides necessary and sufficient conditions for the existence of an extension $g$ to have $c(g)=b(f)$. 
\begin{theorem}\citep{Liu20173}
\label{iff}
Let $(X,\mathcal{I})$ be a matroid and $f: \mathcal{I}\rightarrow \mathbb{R}$  a polymatroid function. Let $g:2^X\rightarrow\mathbb{R}$ be a polymatroid function that agrees with $f$ on $\mathcal{I}$. Then $c(g)=b(f)$ if and only if 
\begin{equation}
\label{sufneccond2}
g(X)-g(X\setminus \{a\})\geq (1-b(f))g(\{a\})
\end{equation}
 for any $a\in  X$, and equality holds for some $a\in X$.
\end{theorem}

\cite{Liu20173} provided the following improved bounds for the greedy strategy if there exists an extension $g$ of $f$ such that $c(g)=b(f)$.
\begin{theorem}\citep{Liu20173}
\label{ImprovedBoundsCor}
Let $(X,\mathcal{I})$ be a matroid of rank $K$. Let $g:2^X\rightarrow\mathbb{R}$ be a polymatroid function that agrees with $f$ on $\mathcal{I}$ such that $g(X)-g(X\setminus \{a\})\geq (1-b(f))g(\{a\})$ for any $a\in X$ with equality holding for some $a\in X$. Then, any greedy solution $G_K$ to problem~$(\ref{setproblem})$ satisfies
 \begin{equation}
 \label{bound12}
 \frac{f(G_K)}{f(O)}\geq \frac{1}{1+b(f)}.
 \end{equation}
 In particular, when  $(X,\mathcal{I})$ is a uniform matroid, any greedy solution $G_K$ to problem~$(\ref{setproblem})$ satisfies
\begin{align}
\label{bound22}
\frac{f(G_K)}{f(O)}&\geq\frac{1}{b(f)}\left[1-\left(1-\frac{b(f)}{K}\right)^K\right]> \frac{1}{b(f)}\left(1-e^{-b(f)}\right).
\end{align}
\end{theorem}
\begin{remark}
The bounds ${1}/({1+b(f)})$ and $(1-\left(1-{b(f)}/{K}\right)^K)/b(f)$ do not depend on sets outside the matroid, so they apply to problems where the objective function is only defined on the matroid, provided that an extension that satisfies the assumptions in Theorem~\ref{iff} exists. When $f$ is defined on the entire power set, $b(f)\leq c(f)$, which implies that the bounds are stronger than those  from \cite{conforti1984submodular}. 
\end{remark}

 Next consider again the task assignment problem from Section~\ref{examples}. \cite{Liu20173} gave an extension $g$ of $f$ defined on the uniform matroid of rank $2$ to the whole power set with $c(g)=b(f)$, which is reviewed as follows.
 
\textbf{Example:} Let $X=\{a_1,a_2,a_3,a_4\}$, $p(a_1)=0.4$, 
$p(a_2)=0.6$, $p(a_3)=0.8$, and $p(a_4)=0.9$. Then, $f(A)$ is defined as in (\ref{objectivefunction}) for any $A=\{a_i,\ldots,a_k\}\subseteq X$. Let $K=2$, then  $\mathcal{I}=\{S\subseteq X: |S|\leq 2\}$. It is easy to show that $f:\mathcal{I}\rightarrow \mathbb{R}$ is a polymatroid function.

The polymatroid function $g$  constructed using (\ref{construction}) while satisfying (\ref{dsatisfy}) and (\ref{sufneccond2}) from \cite{Liu20173} is of the following form:

$g(\{a_1,a_2,a_3\})=f(\{a_2,a_3\})+d_{\{a_1,a_2,a_3\}}=0.96.$

$ g(\{a_1,a_2,a_4\})=f(\{a_2,a_4\})+d_{\{a_1,a_2,a_4\}}=1,$

 $g(\{a_1,a_3,a_4\})=f(\{a_3,a_4\})+d_{\{a_1,a_3,a_4\}}=1.02,$
 
 $ g(\{a_2,a_3,a_4\})=f(\{a_3,a_4\})+d_{\{a_2,a_3,a_4\}}=1.04,$
 
 $g(X)=g(\{a_2,a_3,a_4\})+d_X=1.08$.

The total curvature $c$ of $g:2^X\rightarrow\mathbb{R}$ is
\begin{align*}
c(g_2)=&\max_{a_i\in X}\left\{1-\frac{g(X)-g(X\setminus\{a_i\})}{g(\{a_i\})-g(\emptyset)}\right\}=0.9=b(f)<c(f)=0.992.
\end{align*}

By Theorem~\ref{ImprovedBoundsCor},  the greedy strategy for the task scheduling problem satisfies the bound $(1-(1-{b(f)}/{2})^2)/b(f)=0.775$, which is better than the previous bound 
$(1-(1-{c(f)}/{2})^2)/c(f)=0.752$.

\subsection{Batch Actions}
\label{batchactions}
Suppose we batch the selected actions into batches of size $k$. What results is the \emph{$k$-batch greedy strategy}, which starts with the empty set and iteratively adds to the current solution set a batch of elements with the largest gain in the objective function under the constraints. The greedy strategy we considered in Sections~\ref{greedybounds}--\ref{curvature} is a special case of the batched greedy with batch size equal to $1$. Intuitively, larger $k$ should result in better performance, albeit at the expense of increasing computational complexity. But how do the previous bounds improve as a function of $k$? In this section, we review performance bounds for the $k$-batch greedy strategy.

We start by introducing the $k$-batch greedy strategy  as follows. Consider again problem~(\ref{setproblem}) and write the maximal cardinality of the  sets in $\mathcal{I}$ as $K=k(l-1)+m$, where $l,m$ are nonnegative integers and $0<m\leq k$. Note that $m$ is not necessarily the remainder of $K/k$, because $m$ could be equal to $k$. This happens when $k$ divides $K$. The $k$-batch greedy strategy is as follows \citep{Liu2016, Liu2018}:

\noindent Step~1: Let $S^0=\emptyset$ and $t=0$.

\noindent Step~2: Select $J_{t+1}\subseteq X\setminus S^t$ such that $|J_{t+1}|=k$, $S^t\cup J_{t+1}\in\mathcal{I}$, and 
\begin{align*}
f(S^t\cup J_{t+1})=\max\limits_{J\subseteq X\setminus S^t\ \text{and}\ |J|=k }f(S^t\cup J);
\end{align*}
 then set $S^{t+1}=S^t\cup J_{t+1}$.

\noindent Step~3:  If $t+1< l-1$, set $t=t+1$, and repeat Step~2.

\noindent Step~4: If $t+1=l-1$, select $J_{l}\subseteq X\setminus S^{l-1}$ such that $|J_{l}|=m$, $S^{l-1}\cup J_{l}\in\mathcal{I}$, and 
\[f(S^{l-1}\cup J_{l})=\max\limits_{J\subseteq X\setminus S^{l-1}\ \text{and}\ |J|=m}f(S^{l-1}\cup J).\]

\noindent Step~5: Return the set $S=S^{l-1}\cup J_{l}$ and terminate.

Any set generated by the above procedure is called a \emph{$k$-batch greedy solution}. For the above strategy, there are $l$ steps in total, and exactly $k$ actions are selected at each of the first $l-1$ steps but the final step may select fewer than $k$ actions.   A similar batched greedy strategy is investigated by \cite{hausmann1980} called the \emph{$(\leq k)$-greedy strategy}, where at most $k$ actions are selected at each stage.

The performance of the $k$-batch greedy strategy under uniform matroid constraints was first investigated by \cite{Nemhauser19781}, stated as follows.

\begin{theorem}\citep{Nemhauser19781}
If $(X,\mathcal{I})$ is a uniform matroid of rank $K$ and $f$ is a polymatroid set function, then any $k$-batch greedy solution $S$ satisfies 
\[\frac{f(S)}{f(O)}\geq 1-\left(1-\frac{m}{kl}\right)\left(1-\frac{1}{l}\right)^{l-1}.\]
\end{theorem}
\begin{remark}
When $m=k$, i.e., the batch size $k$ divides the rank $K$, the bound is tight; see \cite{Nemhauser19781} for proof.
\end{remark}
By introducing the \emph{total $k$-batch curvature} 
\begin{equation}
\label{totalkbatchcurvature}
c_k:=\max\limits_{I\in \hat{X}}\left\{1-\frac{\varrho_I(X\setminus I)}{\varrho_I(\emptyset)}\right\},
\end{equation}
 where $\hat{X}=\{I\subseteq X: \varrho_I(\emptyset)\neq 0 \ \text{and}\ |I|=k\}$,
 \cite{Liu2018} derived performance bounds in terms of $c_k$ for the $k$-batch greedy strategy under both general matroid and uniform matroid constraints, and investigated the monotoneity of the performance bounds with respect to the batch size $k$.
\begin{theorem}\citep{Liu2018}
Assume that $f$ is a polymatroid set function. When $(X,\mathcal{I})$ is a general matroid, then any $k$-batch greedy solution satisfies 
\[\frac{f(S)}{f(O)}\geq \frac{1}{1+c_k}.\]
When $(X,\mathcal{I})$ is a uniform matroid, then any $k$-batch greedy solution satisfies 
\[\frac{f(S)}{f(O)}\geq \frac{1}{c_k}\left[1-\left(1-\frac{c_k}{l}\frac{m}{k}\right)\left(1-\frac{c_k}{l}\right)^{l-1}\right].\]

\end{theorem}
\begin{remark}
When $k=1$, the bound for general matroid constraint becomes $1/(1+c)$ and the bound for uniform matroid constraints becomes $(1-(1-c/K)^K)/c$, which is consistent with the results in Theorem~\ref{greedyboundcurvature}.
\end{remark}
\begin{remark}
\label{curvaturemonotone}
The total $k$-batch curvature is nonincreasing in $k$, i.e.,  $c_{k_2}\leq c_{k_1}$ whenever $k_2\geq k_1$ \citep{Liu2018}.
\end{remark}
\begin{remark}
Based on Remark~\ref{curvaturemonotone}, we can discuss the monotoneity of the bounds for both general matroid and uniform matroid constraints. The bound $1/(1+c_k)$ for general matroid constraints is monotone in $k$. For uniform matroid constraints, when the batch size $k$ divides $K$, the bound becomes \[\frac{1}{c_k}\left[1-\left(1-\frac{c_k}{l}\right)^{l}\right],\]
which is monotone in $k$. Moreover, 
\[\frac{1}{c_k}\left[1-\left(1-\frac{c_k}{l}\right)^{l}\right]> \frac{1-e^{-c_k}}{c_k}\geq \frac{1}{1+c_k},\]
which means that the bound for uniform matroid constraints is better than the bound for general matroid constraints. However, if $k$ does not divide $K$, the exponential bound might be worse than the harmonic bound. For example, when $K=100, k=80$, and $c_k=0.6$,  the exponential bound is $0.5875$, which is worse than the harmonic bound  $0.6250$ \citep{Liu2018}.
\end{remark}

\textbf{Examples:} Now consider again the task assignment and adaptive sensing problems from Section~\ref{examples} to demonstrate that the total curvature $c_k$ decreases in $k$ and the performance bound for a uniform matroid  increases in $k$ under the condition that the batch size $k$ divides the rank $K$.

\textbf{Task Assignment Problem}:
We still  order the elements of $X$ as $a_{[1]}, a_{[2]}, \ldots, a_{[N]}$ such that $$0<p(a_{[1]})\leq p(a_{[2]})\leq\ldots\leq p(a_{[N]})\leq 1.$$ Then by the definition of the total curvature $c_k$, we have 
\[
c_k=\max\limits_{i_1,\ldots,i_k\in {X}}\left\{1-\frac{f(X)-f(X\setminus\{i_1,\ldots,i_k\})}{f(\{i_1,\ldots, i_k\})-f(\emptyset)}\right\}=1-\prod_{l=k+1}^N(1-p(a_{[l]})).
\]
From the expression of $c_k$, we can see that $c_k$ is nonincreasing in $k$, but when $N$ is large, $c_k$ is close to 1 for each $k$.

To numerically evaluate the relevant quantities here, \cite{Liu2018} randomly generated a set  $\{p(a_i)\}_{i=1}^{30}$. In Fig.~3, they considered $K=20$,  and  batch sizes $k=1,2,\ldots, 10$.   Fig.~3 shows that  the exponential bound for $k=3,6,8,9$ is worse than that for $k=1,2$, which illustrates our earlier remark that the exponential bound for the uniform matroid case is not necessarily 
monotone in $k$  even though $c_k$ is monotone in $k$. Fig.~3 also shows that  the exponential bound $(1-(1-c_k/l\cdot m/k)(1-c_k/l)^{l-1}/c_k$ coincides with $(1-(1-c_k/l)^{l})/c_k$ for $k=1,2,4,5,10$ and it is nondecreasing in $k$, which illustrates our remark that  the exponential bound is nondecreasing in $k$ under the condition that $k$ divides $K$. Owing to the nature of the total curvature for this example, it is not easy to see that $c_k$ is nonincreasing in $k$ (all $c_k$ values here are very close to 1). 
\begin{figure}
   \centering
   \includegraphics[width=0.65\columnwidth, trim= 30 120 0 130, clip]{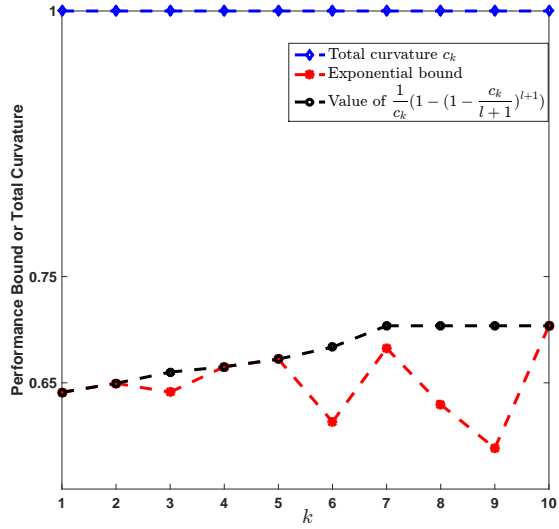}
   \caption{Total curvature/performance bounds for greedy strategy in task assignment problem}
\label{Fig3}
\end{figure}

\textbf{Adaptive Sensing}: For convenience, set $\sigma=1$.
Then, we have 
\begin{align*}
c_k&=\max\limits_{J_k\subseteq X, |J_k|=k}\left\{1-\frac{f(X)-f(X\setminus J_k)}{f(J_k)}\right\}\nonumber\\
&=\max\limits_{J_k\subseteq X, |J_k|=k}\left\{1- \frac{\log (st)-\log\left(s-\sum\limits_{i: e_i\in J_k}e_i\right)\left(t-\sum\limits_{i: e_i\in J_k}(1-e_i)\right)}{\log\left(1+\sum\limits_{i: e_i\in J_k}e_i\right)\left(1+\sum\limits_{i: e_i\in J_k}(1-e_i)\right)}\right\},
\end{align*}
where $X=\{B_1,\ldots, B_N\}$, $s=1+\sum_{i=1}^Ne_i$, and $t=1+\sum_{i=1}^N(1-e_i)$.

We already saw that    the exponential bound for the uniform matroid case is not necessarily monotone in $k$ from the task assignment problem, so we will only consider the case when the batch size $k$ divides $K$. \cite{Liu2018} considered  $ K=24$ for $k=1,2,3,4,6,8$ in Fig.~4. The figure shows that the curvature decreases in $k$ and the exponential bound increases in $k$ since $k$ divides $K$ for $k=1,2,3,4,6,8$, which again demonstrates the claim that $c_k$ decreases in $k$ and the exponential bound increases in $k$ under the condition that $k$ divides $K$. 
\begin{figure}
   \centering
   \includegraphics[width=\columnwidth, trim= 0 0 0 0, clip]{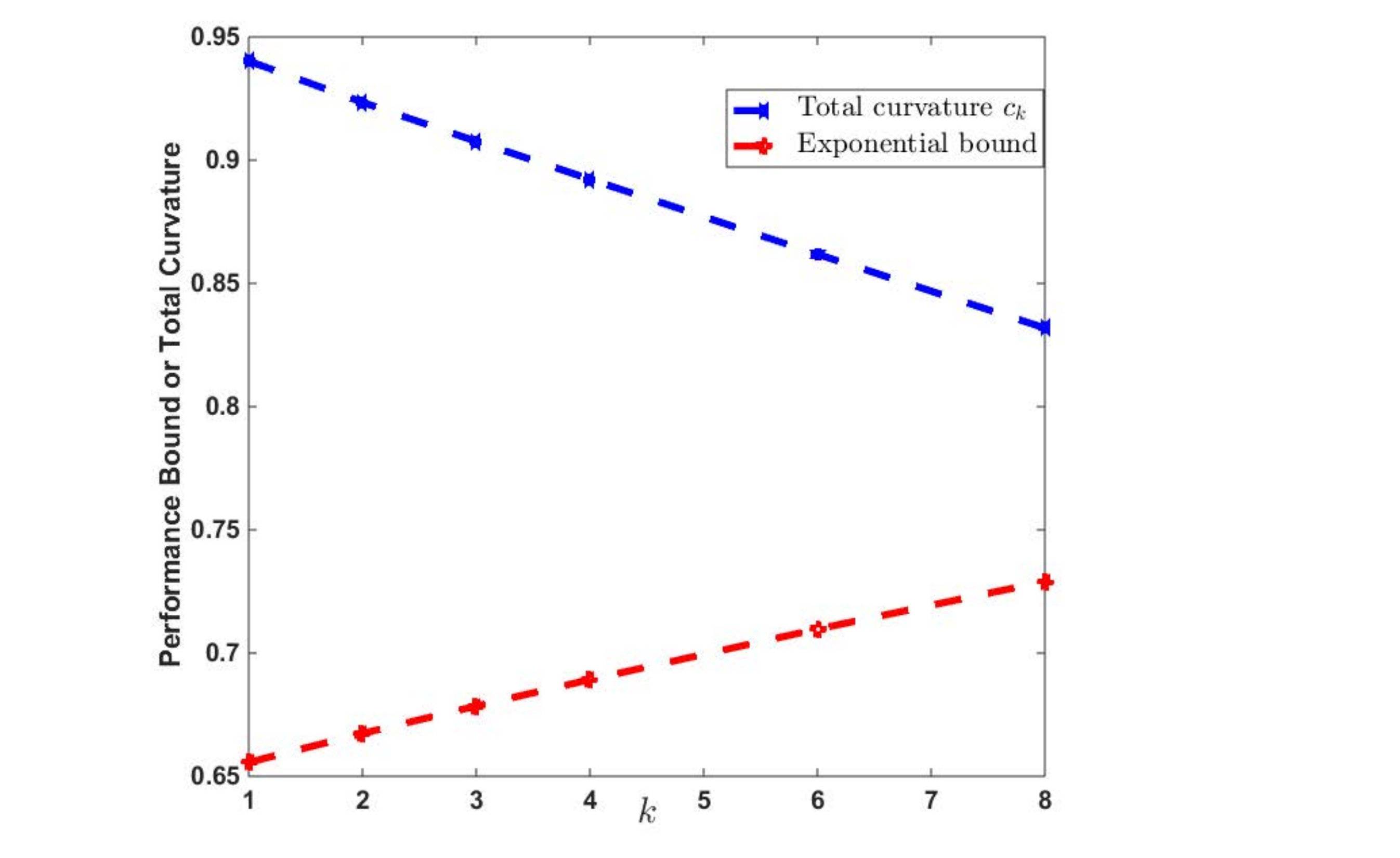}
   \caption{Total curvature/performance bound for greedy strategy in adaptive sensing problem}
\label{Fig4}
\end{figure}

\subsection{Noncooperative Games}
\label{games}
In the previous sections, we reviewed performance bounds for greedy-type strategies in set submodular optimization problems. It turns out that similar techniques can be used to bound the performance of Nash equilibria in noncooperative games--utility maximization problems. The connection to the game setting is easy to imagine by associating the objective function in set optimization  with a social utility function in games, greedy strategies with Nash equilibria, and batching with cooperation of subgroups in games. We first introduce some background on utility maximization problems and Nash equilibria.

A great number of interesting practical problems can be posed as utility maximization problems: these include facility location  \citep{Atamturk2011}, traffic routing and congestion management  \citep{Marden2007,Rexford2007}, sensor selection  \citep{Portal2007,Liu2014}, and network resource allocation  \citep{Anantharam2002,Chiang2007}. In a utility maximization problem, a set of users make decisions according to their own set of feasible strategies, resulting in an overall social utility value, such as profit, coverage, achieved data rate, and quality of service. The goal is to maximize the social utility function. Often, the users do not cooperate in selecting their strategies.

In general, it is impractical to find the  optimal strategy maximizing the social utility function. However, it is feasible to consider scenarios where individual users or groups of users separately maximize their own private objective functions. The usual framework for studying such scenarios is game theory together with its celebrated notion of Nash equilibria.
A \emph{Nash equilibrium} is a set of strategies (deterministic or randomized) for which no user can improve its own private utility by changing its strategy unilaterally. \cite{Nash1951} proved that any finite and non-cooperative game has at least one Nash equilibrium. 

The performance of Nash equilibria compared with the  optimal solution in submodular utility maximization problems was investigated by \cite{Vetta2002}.  Based on the existing results, \cite{Liu20172} established bounds for Nash equilibria when there is ``grouping" among users, which is useful in understanding the role of cooperation and social ties in games. Before we review these results, we  introduce some notation and terminology from \cite{Vetta2002} and \cite{Liu20172}.

Suppose we have a set  $\mathcal{N}=\{1,2,\ldots,N\}$ of $N$ users. Each element in $V_i$ $(i=1,\ldots, N)$ represents an \emph{act} that user $i$ can take. We call a set of acts an \emph{action}, and if an action $x_i\subseteq V_i$ is available to user $i$ we call it a \emph{feasible action}. We denote by $\mathcal{X}_i$ the set of all feasible actions for user $i$, i.e., $\mathcal{X}_i=\{x_i\subseteq V_i: x_i$ is a feasible action$\}$, with $n_i=|\mathcal{X}_i|$ the cardinality of $\mathcal{X}_i$. We call $\mathcal{X}_i$ the \emph{action space} for user $i$.  A \emph{pure strategy} is one in which the user takes a specific
action. A \emph{mixed strategy} is one in which the user takes actions according to some probability distribution. The set of mixed strategies is called the \emph{strategy space}. We represent the strategy space for user $i$ by $\mathcal{S}_i=\{s_i\in\mathbb{R}^{n_i}:\sum_{j=1}^{n_i}s_i^j=1,s_i^j\geq 0\}$, where $s_i=(s_i^1,\ldots, s_i^{n_i})$ is called a \emph{strategy taken by user $i$} and $s_i^j\ge 0$ is the probability with which user $i$ takes action $j$. When $s_i^j=1$ for some $j$ and $s_i^l=0$ for all $l\neq j$, user $i$ is said to take a \emph{pure} strategy. Otherwise,  user $i$ takes a \emph{mixed} strategy. Write $\mathcal{S}=\prod_{i=1}^N\mathcal{S}_i$. The indexed set $S=(s_{1},\ldots, s_{N})$, with $s_{i}\in\mathcal{S}_i$ and $i=1,\ldots,N$, is called a \emph{strategy set} of size $N$ in $\mathcal{S}$. 

Given a strategy set $S=(s_1,\ldots,s_N)\in\mathcal{S}$, the set $S_{-i}=(s_1,\ldots,s_{i-1},s_{i+1},\ldots,s_N)$ is the subset of $S$ that contains strategies taken by all users except user $i$, and $(S_{-i},s_{i}')=(s_1,\ldots,s_{i-1},s_i',s_{i+1},\ldots,s_N)$ is the strategy set that results from $S$ when user $i$ changes its strategy from $s_i$ to $s_i'$.  %Given strategy sets $T=(t_{i_1},\ldots, t_{i_k})$  and $W=(w_{j_1},\ldots, w_{j_l})$,  $T\oplus W=(t_{i_1},\ldots, t_{i_k},w_{j_1},\ldots, w_{j_l})$ denotes the concatenation of $T$ and $W$ when $i_p\neq j_q$ for $1\leq p\leq k$ and $1\leq q \leq l$.

The \emph{expected  social utility} function and \emph{expected  private utility} function for user $i$  from strategies in $\mathcal{S}$ to real numbers are denoted by $\bar{\gamma}$ and $\bar{\alpha}_i$, respectively. Define $\bar{\gamma}_{s_i}(S_{-i})=\bar{\gamma}(S)-\bar{\gamma}(S_{-i})$ for any set $S=(s_1,\ldots, s_N)\in \mathcal{S}$ and $s_i$ $(i=1,\ldots, N)$.

Now we introduce the definition of a Nash equilibrium and a valid system, then review performance bounds for Nash equilibria under some conditions from \cite{Vetta2002}.
\begin{definition}
\label{dfn:NE}
A strategy set $S\in\mathcal{S}$ is a \emph{Nash equilibrium} if no user has an incentive to unilaterally change its strategy, i.e., for any user $i$, 
\begin{equation}
\label{ineq:NE}
\bar{\alpha}_i(S)\geq \bar{\alpha}_i((S_{-i}, s_i')),\quad \forall s_i'\in\mathcal{S}_i.
\end{equation}
\end{definition}

\medskip
\begin{assumption}
\label{assumption1}
\citep{Vetta2002} The private utility of user $i$ ($i=1,\ldots,N$) is at least as large as the loss in the social utility resulting from user $i$ dropping out of the game. That is, the system ($\bar{\gamma}, \{\bar{\alpha}_i\}_{i=1}^N$) has the property  that  for any strategy set $S=(s_1,\ldots, s_N)\in\mathcal{S}$,
\begin{equation}
\label{assum1}
\bar{\alpha}_i(S)\geq \bar{\gamma}_{s_i}(S_{-i}),\quad \forall i=1,\ldots,N.
\end{equation}
\end{assumption}

\medskip
\begin{assumption}
\label{assumption2}
\citep{Vetta2002} The sum of the private utilities of the system is not larger than the social utility, i.e., for any strategy set $S=(s_1,\ldots, s_N)\in\mathcal{S}$,
\begin{equation}
\label{assum2}
\sum\limits_{i=1}^N\bar{\alpha}_i(S)\leq\bar{\gamma}(S) .
\end{equation}
\end{assumption}

A utility system $(\bar{\gamma}, \{\bar{\alpha}_i\}_{i=1}^N)$ satisfying Assumptions~\ref{assumption1} and~\ref{assumption2} is called a \emph{valid} system. 
We denote by $\Omega=(\omega_1,\ldots,\omega_N)$ the optimal strategy set in maximizing an expected utility function $\bar{\gamma}$, and assume that $\Omega$ is composed of pure strategies $\omega_i\in\mathcal{S}_i$, $i=1,\ldots, N$. For convenience, we also use $\omega_i$ to denote the optimal action that user $i$ takes.  Consider a strategy set $S=(s_1,\ldots,s_i)$ where $i=1,\ldots,N$. Suppose that user $j$ ($j=1,\ldots,i$) uses a mixed strategy $s_j$ that takes actions $x_j^1,\ldots, x_j^{n_j}$ with probabilities $s_j^1,\ldots, s_j^{n_j}$.
We use the notation $\Omega\cup S$ to represent the strategy in which user $j$ $(j=1,\ldots, i)$ takes the actions $\omega_j\cup x_j^1,\ldots,\omega_j\cup x_j^{n_j}$ with probabilities $s_j^1,\ldots,s_j^{n_j}$, and user $j$ $(j=i+1,\ldots, N)$ takes the action $\omega_j$, so $\Omega\cup S$ is well defined.
\begin{theorem}
\citep{Vetta2002} For a valid utility system $(\bar{\gamma}, \{\bar{\alpha}_i\}_{i=1}^N)$, if the expected social utility function $\bar{\gamma}$ is submodular, then for any Nash equilibrium $S\in\mathcal{S}$ we have 
\begin{equation}
\label{ineq:NE}
\bar{\gamma}(S)\geq\frac{1}{2}\left(\bar{\gamma}(\Omega)+\sum_{i=1}^N\bar{\gamma}_{s_i}({\Omega\cup S_{-i}})\right).
\end{equation}
\end{theorem}
\begin{remark}
If $\bar{\gamma}$ is monotone, then $\bar{\gamma}_{s_i}({\Omega\cup S_{-i}})\geq 0$ and the above inequality shows that any Nash equilibrium achieves at least $1/2$ of the optimal social utility function value.
\end{remark}
By defining the curvature $c$ of the expected social utility function $\bar{\gamma}$,
\[c:=\max_{i:\bar{\gamma}_{s_i}(\emptyset)\neq 0}\left\{1-\frac{\bar{\gamma}_{s_i}(\Omega\cup S_{-i})}{\bar{\gamma}_{s_i}(\emptyset)}\right\},\] 
\cite{Vetta2002} derived the following tighter performance bound in terms of the curvature for  Nash equilibria.
\begin{theorem}
\label{thm:NEC}
\citep{Vetta2002} For a valid utility system $(\bar{\gamma}, \{\bar{\alpha}_i\}_{i=1}^N)$, if  the expected social utility function $\bar{\gamma}$ is monotone and submodular, then for any Nash equilibrium $S\in\mathcal{S}$ we have 
\begin{equation}
\label{ineq:NEC}
\frac{\bar{\gamma}(S)}{\bar{\gamma}(\Omega)}\geq\frac{1}{1+c}.
\end{equation}
\end{theorem}
\begin{remark}
When the expected social utility function $\bar{\gamma}$ is monotone and submodular, we have $c\in[0,1]$, which implies that $
\bar{\gamma}(S)\geq \bar{\gamma}(\Omega)/2$.
\end{remark}

Next we review performance bounds for \emph{group Nash equilibria} 
defined by \cite{Liu20172}. They considered the case where the set of all users in the utility maximization system are divided into disjoint groups, and the users in the same group choose their strategies by maximizing their group utility function jointly.

Assume that the set of users $\mathcal{N}=\{1,\ldots, N\}$ is divided into $l$ disjoint groups, in which group $i$ ($i=1,\ldots,l$) has users $\{m_i+1,\ldots,m_i+k_i\}$, where $m_i=\sum_{j=1}^{i-1}k_j$, $k_j$ is the number of users in group $j$, and $\sum_{j=1}^lk_j=N$.  Let $s^i=(s_{m_i+1},\ldots, s_{m_i+k_i})$ denote the  \emph{group strategy} for group $i$, where 
$s_i\in\mathcal{S}_i$ is the strategy for user $i$.  This includes the strategies taken by all the users in group $i$ ($i=1,\ldots,l$). Let $S^{-i}$  denote the set of group strategies taken by all groups except for group $i$ and $(S^{-i}, t^i)$ denote the group strategy set obtained when group $i$ changes its group strategy from $s^i$ to $t^i$. Let $\bar{\eta}_i$ denote the expected group utility function for group $i$. Define $\bar{\gamma}_{s^i}(S^{-i})=\bar{\gamma}(S)-\bar{\gamma}(S^{-i})$ for any $S=(s^1,\ldots, s^l)\in\mathcal{S}$ and $s^i$ ($i=1,\ldots,l$).

\begin{definition}
 A strategy set $S=(s^1,\ldots, s^l)\in\mathcal{S}$ is a \emph{group Nash equilibrium} of a utility system if no group  can improve its group utility by unilaterally changing its group strategy, i.e., for any $i=1,\ldots,l$,
\[
 \bar{\eta}_i(S)\geq \bar{\eta}_i((S^{-i}, t^{i})), \quad \forall t^i=(t_{m_i+1},\ldots, t_{m_i+k_i}),
\]
where $t_j\in\mathcal{S}_j$ for $j=m_i+1,\ldots, m_i+k_i$. 
 
 \end{definition}

The utility system $(\bar{\gamma}, \{\bar{\eta}_i\}_{i=1}^l)$ is \emph{valid} if it satisfies the following two assumptions \citep{Liu20172}. 
 
\begin{assumption}
\label{assumption5}
The group utility of group $i$ is at least as large as the loss in the social utility resulting from all the users in group $i$ dropping out of the game. That is, the system $(\bar{\gamma}, \{\bar{\eta}_i\}_{i=1}^l)$  has the property that for any strategy set $S=(s^1,\ldots, s^l)\in\mathcal{S}$,
\begin{equation}
\label{ineq:ass5}
\bar{\eta}_i(S)\geq\bar{\gamma}_{s^i}(S^{-i}),\quad \forall i=1,\ldots,l.
\end{equation}
\end{assumption}

\medskip
\begin{assumption}
\label{assumption6}
The sum of the group utilities of the system is not larger than  the social utility, i.e., for any strategy set $S=(s^1,\ldots, s^l)\in\mathcal{S}$,
\begin{equation}
\label{ineq:ass6}
\sum\limits_{i=1}^l\bar{\eta}_i(S)\leq \bar{\gamma}(S).
\end{equation}
\end{assumption}

\medskip
%A utility system $(\gamma, \{\eta_i\}_{i=1}^N)$  is a valid system if it satisfies Assumptions~\ref{assumption5} and \ref{assumption6}.

\begin{theorem}\citep{Liu20172}
\label{thm:GNE}
For a valid utility system $(\bar{\gamma}, \{\bar{\eta}_i\}_{i=1}^l)$, if the expected social utility function $\bar{\gamma}$ is submodular, then  any group Nash equilibrium $S=(s^1,\ldots,s^l)\in\mathcal{S}$ satisfies
\begin{equation}
\label{ineq:GNS}
\bar{\gamma}(S)\geq \frac{1}{2}\left(\bar{\gamma}(\Omega)+\sum\limits_{i=1}^l\bar{\gamma}_{s^i}(\Omega\cup S^{-i})\right).
\end{equation}
\end{theorem}
To better characterize the relation of the social utility value of any group Nash equilibrium and that of the optimal solution $\Omega$, \cite{Liu20172} defined the {group curvature} $c_{k_i}$ of the social utility function for group $i$ as 
\[c_{k_i}:=\max\limits_{S\in\mathcal{S}, \bar{\gamma}_{s^i}(\emptyset)\neq 0}\left\{1-\frac{\bar{\gamma}_{s^i}(\Omega\cup S^{-i})}{\bar{\gamma}_{s^i}(\emptyset)}\right\}{.}\]
\begin{theorem}\citep{Liu20172}
\label{thm:GNEC}
For a valid utility system $(\bar{\gamma}, \{\bar{\eta}_i\}_{i=1}^l)$, if the expected social utility function $\bar{\gamma}$ is monotone and submodular, then  any group Nash equilibrium $S=(s^1,\ldots,s^l)\in\mathcal{S}$ satisfies
\[\bar{\gamma}(S)\geq  \frac{1}{1+\max\limits_{1\leq i\leq l}c_{k_i}}\bar{\gamma}(\Omega).\]
In particular, if $\mathcal{X}_1=\mathcal{X}_2=\cdots=\mathcal{X}_N$, we have 
\[\bar{\gamma}(S)\geq  \frac{1}{1+c_{k^*}}\bar{\gamma}(\Omega),\]
where $k^*=\min_{1\leq i\leq l}k_i$.
\end{theorem}

\begin{remark}
 When the expected group utility function $\bar{\gamma}$ is monotone and submodular, it is easy to check that $c_{k_i}\in[0,1]$, which implies that $1/(1+\max_{1\leq i\leq l}c_{k_i})\geq 1/2$.
\end{remark}
\begin{remark}
  When the expected group utility function $\bar{\gamma}$ is monotone and submodular, we have $\bar{\gamma}(S)\geq\bar{\gamma}(\Omega)/ (1+\max_{1\leq i\leq l}c_{k_i})\geq \bar{\gamma}(\Omega)/ (1+c)$. This shows that the bound for the case with grouping is tighter than that for the case without grouping. Of course, this is unsurprising, because grouping entails cooperation. Moreover, under the condition that each user has the same strategy space,  the larger the value of $k_i$, the higher the degree of cooperation, and the tighter the lower bound. 
\end{remark}

\section{Strings of Actions}
\label{stringaction}
In Section~\ref{setactions}, we considered the optimization problem where the argument of  the objective function is a set of actions. Suppose the objective function depends not only on the set of actions but also on the order of actions. We call the argument of the objective function a \emph{string of actions.} In this section, we introduce notation and terminology for strings and string functions, string optimization problem, performance bounds for the greedy strategy, and applications.
\subsection{Notation and Terminology}
\label{stringnotation}
Let $X$ be a set of all possible actions. We use $A=(a_1,a_2,\ldots,a_k)$ ($a_i\in X$) to denote a \emph{string} of actions taken over $k$ consecutive stages. We define its \emph{length} as $k$, denoted by $|A|=k$. Note that $k=0$ corresponds to the empty string, denoted by $A=\emptyset$.

Let ${X}^*$ denote the set of all possible strings of actions. If two strings in ${X}^*$ are expressed by $M=(a_1^m,a_2^m,\ldots, a_{k_1}^m)$ and $N=(a_1^n,a_2^n,\ldots, a_{k_2}^n)$, we write $M=N$ iff $k_1=k_2$ and $a_i^m=a_i^n$ for each $i=1,2,\ldots,k_1$. Moreover, we define string {concatenation} as $M\oplus N= (a_1^m,a_2^m,\ldots, a_{k_1}^m,a_1^n,a_2^n,\ldots, a_{k_2}^n)$. 

We write $M\preceq N$ if we have $N=M\oplus L$ for some $L\in X^*$. In this case, we also say that $M$ is a \emph{prefix} of $N$.
We write $M\prec N$ if there exists a set of strings $L_i\in X^*$ such that 
$N=L_1\oplus(a_1^m,\ldots, a_{i_1}^m)\oplus L_2 \oplus (a_{i_1+1}^m,\ldots, a_{i_2}^m)\oplus\cdots\oplus(a_{i_{k-1}+1}^m,\ldots, a_{k_1}^m)\oplus L_k$. Note that $\prec$ is weaker than $\preceq$, which means $M\preceq N$ implies $M\prec N$, but the converse is not necessarily true.
%
%Assume that there exists $K$ such that for all $M\in\mathcal{I}$ we have $|M|\leq K$ and there exists a $N\in\mathcal{I}$ such that $|N|=K$. The string $N$ is called a \emph{maximal string}.
%
%For a given $K$, let $\mathcal{I}=\{A\in\mathbb{A}^*: |A|\leq K\}$, the pair $(\mathbb{A},\mathcal{I})$ is called a  \emph{uniform string-matroid}.

Similar to the definition of a polymatroid set function in Section~\ref{setfunctions}, we define a function from strings to real numbers, $f: X^*\to \bbR$,  a \emph{polymatroid string function}  if
\begin{itemize}
\item [i.] $f(\emptyset)=0$.
\item[ii.] $f$ has the \emph{prefix-monotone} property: $ \forall M, N \in X^*,$  $f(M\oplus N)\geq f(M)$.
\item[iii.] $f$ has the \emph{diminishing-return} property: $\forall M\preceq N \in X^*, \forall a\in X$, $f(M\oplus (a))-f(M) \geq f(N\oplus (a))-f(N)$.
\end{itemize}
A function  $f: X^*\to \bbR$ is \emph{postfix monotone} if
 \[\forall M, N \in X^*, f(M\oplus N)\geq f(N).\]
 Notice the difference between the prefix-monotone property and postfix-monotone property. 

Let $\mathcal{I}$ denote a collection of strings from $X^*$.
The pair $(X,\mathcal{I})$ is called a \emph{string matroid}  \citep{Zhang2016} if $\mathcal{I}$ satisfies the following properties:
\begin{itemize}
\item[i.] $\mathcal{I}$ is non-empty;

\item[ii.] \emph{Hereditary}: $\forall M\in\mathcal{I}, N\prec M$ implies that $n\in\mathcal{I}$;

\item[iii.] \emph{Augmentation}: $\forall M, N\in\mathcal{I}$ and $|M|<|N|$, there exists an element $x\in X$ in the string $N$ such that $M\oplus (x)\in\mathcal{I}$.
\end{itemize}
The length of the longest string in $\mathcal{I}$ is called the \emph{rank} of $(X,\mathcal{I})$. When $\mathcal{I}=\{A\in X^*: |A|\leq K\}$, the pair $(X,\mathcal{I})$ is called a \emph{uniform string matroid} of rank $K$.

\subsection{String Optimization Problem}
\label{stringoptimization}
In this section, we first formulate the string optimization problem and define the greedy strategy. Then we review performance bounds for the greedy strategy under  uniform string matroid constraints and general string matroid constraints. 

In a variety of  problems in engineering and applied science such as sequential decision making \citep{Littman1996,Roijers2013},  adaptive sensing \citep{Liu2014,Krause2008}, and adaptive control \citep{Jarvis1975,Schlegel2005}, we are faced with optimally choosing a string (ordered set) of actions over a finite horizon to maximize an objective function under some constraints. We call this class of optimization problems string optimization. For set optimization problems, the objective function is not influenced by the order of actions. However, for  string optimization problems, the objective function depends on the order of actions. 
 Let $f: X^*\to \bbR$ be an objective function. The goal is to find a string $M$, with the constraint $M\in\mathcal{I}$, to maximize the objective function:
\begin{align}\label{stringproblem}
\begin{array}{l}
\text{maximize} \ \    f(M), \ \quad \text{subject to} \ \  M\in \mathcal{I},
\end{array}
\end{align}
where   $X^*$ denotes the set of all possible strings of actions and  $\mathcal{I}$ is a collection of strings from $X^*$.

The solution to the string optimization problems can be characterized using  dynamic programming via Bellman's principle \citep{Bertsekas2005,Powell2007}. However, dynamic programming suffers from the curse of dimensionality and is therefore impractical for many problems of interest. Hence, we often turn to approximation techniques. One approximation technique is the greedy strategy, which is to find an action at each stage to maximize the step-wise gain in the objective function. The performance for the greedy strategy in string optimization problems has been  investigated by  \protect\cite{streeter2008}, \protect\cite{Zhang2016}, and \protect\cite{Liu2015}. And these specific results will be reviewed  in this section.

Assume that the rank of $(X,\mathcal{I})$ is $K$.
We now define optimal and greedy strategies for problem \eqref{stringproblem} and some related notation.

\noindent\textbf{Optimal String}: 
Any string $O$ is called an \emph{optimal solution} of  Problem (\ref{stringproblem}) if
 \[O\in\mathop{\argmax}\limits_{M\in \mathcal{I}} f(M).\]
If $f$ is prefix monotone, then there exists at least one optimal string of length $K$, denoted by $O_K=(o_1, \ldots, o_K)$.

\newpage
\noindent\textbf{Greedy Algorithm}: 

\noindent\textbf{Input}: A string matroid  $(X,\mathcal{I})$ of rank $K$, a set function $f:X^*\rightarrow \mathbb{R}$, collection $\mathcal{I}$, size $K$

\noindent\textbf{Output}: A string $G_K\in\mathcal{I}$

\noindent$G_0\leftarrow \emptyset$

\noindent For $i=1,\ldots, K$,

\noindent 1. $g_i\leftarrow \argmax\limits_{\substack{a\in X, G_{i-1}\oplus (a)\in\mathcal{I}}}f(G_{i-1}\oplus (a))$

\noindent 2. $G_i\leftarrow G_{i-1}\oplus (g_i)$

Any output of the above algorithm is called a \emph{greedy solution}.
There may exist more than one greedy solution.

\subsection{Performance Bounds for Greedy Strategy}
\label{greedyperformance}
 \cite{streeter2008} first derived performance bounds for the greedy strategy under  uniform string matroid constraints, stated as follows.
\begin{theorem}\citep{streeter2008}
\label{uniformstringmatroidthm1}
Let $(X,\mathcal{I})$ be a uniform string matroid. If $f: X^*\rightarrow\mathbb{R}$ is a polymatroid string function and postfix monotone, then any greedy string $G_K$ satisfies
\begin{equation}
\label{uniformstringbound1}
\frac{f(G_K)}{f(O_K)}\geq 1-\left(1-\frac{1}{K}\right)^K>1-e^{-1}.
\end{equation}
\end{theorem}
\begin{remark}
The same bound holds if $f$ satisfies $f(G_i\oplus O_K)\geq f(O_K)$ for $i=1,\ldots, K-1$, which is weaker than being postfix monotone.
\end{remark}
 \cite{Zhang2016} investigated performance bounds for the greedy strategy under both uniform string matroid and general string matroid constraints by defining the following curvatures.

The \emph{total backward curvature of $f$} is defined as \citep{Zhang2016}
\begin{align}
\label{backwardcurvature}
\sigma:=\max\limits_{\substack{a\in X, M\in X^*\\ f((a))\neq f(\emptyset)}}
\left\{1-\frac{f((a)\oplus M)-f(M)}{f((a))-f(\emptyset)}\right\}.
\end{align}
When $f$ is postfix monotone and diminishing return, we have $0\leq \sigma\leq 1$. The total backward curvature is an upper bound on the second-order difference, over all possible actions $a$ and strings $M$. Next, \cite{Zhang2016} defined the total backward curvature of $f$ with respect to the optimal string $O_K$ by
\begin{align}
\label{backwardcurvaturewithM}
\sigma(O_K):=\max\limits_{\substack{N\in X^*, 0<|N|\leq K\\ f(N)\neq f(\emptyset)}}
\left\{1-\frac{f(N\oplus O_K)-f(O_K)}{f(N)-f(\emptyset)}\right\}.
\end{align}
When $f$ is postfix monotone and string submodular, it is easy to prove that $0\leq \sigma(O)\leq \sigma\leq 1$.

\begin{theorem}\citep{Zhang2016}
\label{uniformstringmatroidthm2}
Let $(X,\mathcal{I})$ be a uniform string matroid of rank $K$. If $f: X^*\rightarrow\mathbb{R}$ is a polymatroid string function, then any greedy string $G_K$ satisfies
\begin{equation}
\label{uniformstringbound2}
\frac{f(G_K)}{f(O_K)}\geq \frac{1}{\sigma(O_K)}\left[1-\left(1-\frac{\sigma(O_K)}{K}\right)^K\right]>\frac{1}{\sigma(O_K)}\left(1-e^{-\sigma(O_K)}\right).
\end{equation}
Moreover, if $f$ is postfix monotone, then any greedy string $G_K$ satisfies
\begin{equation}
\label{uniformstringbound3}
\frac{f(G_K)}{f(O_K)}\geq \frac{1}{\sigma}\left[1-\left(1-\frac{\sigma}{K}\right)^K\right]>\frac{1}{\sigma}\left(1-e^{-\sigma}\right).
\end{equation}
\end{theorem}
\begin{remark}
When $f$ is polymatroid and postfix monotone, we have $0\leq \sigma\leq 1$ by (\ref{backwardcurvature}). So we have $(1-\left(1-{\sigma}/{K}\right)^K)/{\sigma}\geq 1-(1-1/K)^K$ and $(1-e^{-\sigma})/\sigma>1-e^{-1}$, which implies that Theorem~\ref{uniformstringmatroidthm2} provides better bounds than Theorem~\ref{uniformstringmatroidthm1}.
\end{remark}

\begin{theorem}\citep{Zhang2016}
\label{stringmatroidthm}
Let $(X,\mathcal{I})$ be a  string matroid. If $f: X^*\rightarrow\mathbb{R}$ is a polymatroid string function, then any greedy string $G_K$ satisfies
\begin{equation}
\label{stringmatroidbound}
\frac{f(G_K)}{f(O_K)}\geq \frac{1}{1+\sigma(O_K)}.
\end{equation}
Moreover, if $f$ is postfix monotone, then any greedy string $G_K$ satisfies
\begin{equation}
\label{uniformstringbound3}
\frac{f(G_K)}{f(O_K)}\geq \frac{1}{1+\sigma}.
\end{equation}
\end{theorem}

From Theorems~\ref{uniformstringmatroidthm1} and \ref{uniformstringmatroidthm2}, we can see that all the sufficient conditions obtained so far involve strings of length greater than $K$, even though (\ref{stringproblem}) involves only strings up to length $K$. \cite{Liu2015} derived sufficient conditions, which only involve strings of length at most $K$, to have the same bounds hold for uniform string matroid constraints, by defining the following conditions.

A function $f: X^*\to \bbR$ is $K$-\emph{polymatroid} if 
\begin{itemize}
\item [i.] $f(\emptyset)=0$.
\item[ii.] $f$ is {$K$-monotone:} 
$\forall M, N \in X^*,$ and $|M|+|N|\leq K$,  $f(M\oplus N)\geq f(M)$.

\item[iii.] $f$ is $K$-{diminishing:} 
$\forall M\preceq N \in X^*$ and $|N|\leq K-1$, $\forall a\in X$, $f(M\oplus (a))-f(M) \geq f(N\oplus (a))-f(N)$.
\end{itemize}

%Let $G_i=(g_1,\ldots, g_i)$ (as before) and $\bar{O}_{K-i}=(o_{i+1},\ldots, o_K)$ for $i=1,\ldots, K$. Then, $f$ is $K$-\emph{GO-concave} if for $1\leq i\leq K-1$, 

 Let $G_K=(g_1,\ldots, g_K)$ and $\bar{O}_{K-i}=(o_{i+1},\ldots, o_K)$ for $i=1,\ldots, K$. Then, $f$ is $K$-\emph{GO-concave} \citep{Liu2015} if for $1\leq i\leq K-1$, 
\[
f(G_i\oplus \bar{O}_{K-i})\geq \frac{i}{K}f(G_i) + \left(1-\frac{i}{K}\right)f(O_K).
\]
\begin{theorem}\citep{Liu2015}
\label{thm:myopicbounds3}
Let $(X,\mathcal{I})$ be a uniform string matroid. If $f$ is $K$-polymatroid,
then any greedy string satisfies
\[
\frac{f(G_K)}{f(O_K)} \geq \left[1-\left(1-\frac{1}{K}\right)^K\right]>(1-e^{-1}).
\]
\end{theorem}
By defining the curvature $\eta$,
\begin{align*}
\eta=\max\limits_{1\leq i\leq K-1}\left\{\frac{Kf(G_i)-(Kf(G_i\oplus \bar{O}_{K-i})-(K-i)f(O_K))}{(K-i)f(G_i)}\right\},
\end{align*}
 \cite{Liu2015} derived more general performance bounds in
terms of the curvature.
\begin{theorem}\citep{Liu2015}
\label{thm:myopicbounds4}
Let $(X,\mathcal{I})$ be a uniform string matroid. If $f$ is $K$-polymatroid and $K$-GO-concave, then any greedy string satisfies
\begin{align*}
\frac{f(G_K)}{f(O_K)} \geq \frac{1}{\eta}\left[1-\left(1-\frac{\eta}{K}\right)^K\right]>\frac{1}{\eta}(1-e^{-\eta}).
\end{align*}
\end{theorem}
\begin{remark}
If $f$ is $K$-GO-concave, then  we have 
 $0\leq\eta\leq 1$. 
\end{remark}

\textbf{Examples:}
We again consider the task assignment problem and adaptive sensing problem from Section~\ref{examples} to give some sufficient conditions on the parameters of the problems to achieve the performance bound $(1-(1-1/K)^K)$.

\textbf{Task Assignment Problem}:
We use  $p_{i}^j(a)$ to denote the probability of
accomplishing subtask $i$ at stage $j$ when it is assigned to agent
$a\in X$.  Let $a_j$ be the index of the agent selected at stage $j$, the objective function $f$
becomes
\begin{align*}
f((a_1,\ldots, a_k))=\frac{1}{n}\sum\limits_{i=1}^n\left(1-\prod\limits_{j=1}^k\left(1-p_i^j(a_j)\right)\right).
\end{align*} 
For simplicity, we consider the case of $n=1$ (our results can easily be generalized to the case where $n>1$). For $n=1$, the objective function $f$ reduces to   
 \begin{equation}\label{eq:Ex1fn1}
f((a_1,\ldots, a_k))=1-\prod\limits_{j=1}^k(1-p_1^j(a_j)),
\end{equation} 
{and from here on we simply use $p^j(a_j)$ in place of $p_1^j(a_j)$.} 

Note that the value of $f$ depends on the order of the agents selected over time when the probabilities vary from stage to stage. In other words, suppose that we have two agents, Alice and Bob. Then, in general, $p^1(\text{Alice}) \neq p^2(\text{Alice})$, $p^1(\text{Bob}) \neq p^2(\text{Bob})$, $p^1(\text{Alice}) \neq p^1(\text{Bob})$, and $p^2(\text{Alice}) \neq p^2(\text{Bob})$. Therefore, $f((\text{Alice}, \text{Bob})) \neq f((\text{Bob}, \text{Alice}))$.

It is easy to check that $f$ is $K$-monotone and $f(\emptyset)=0$. 

Assume that $p^{j}(a)\in [L(a),U(a)]$, where $L(a) = \min_jp^j(a)$ and $U(a) = \max_jp^j(a)$.
By  \cite{Zhang2016}, a sufficient condition for $f$ to be diminishing return is
\begin{equation}\label{eqn:Ex1Extra}
p^1(g_1)\geq 1-c^K,
\end{equation}
where 
\[
c=\min\limits_{a\in\mathbb{A}}\frac{1-U(a)}{1-L(a)}.
\]
Let $\hat{U} = \max_a{U(a)}$ and $\hat{L}=\min_a{L(a)}$. By \cite{Liu2015}, a sufficient condition for $f$ to be $K$-diminishing is
\begin{equation}\label{eqn:4}
\hat{L}\geq (1-\hat{L})\hat{U},
\end{equation}
and a sufficient condition for $K$-GO-concavity is
\begin{equation}
\label{eqn:Ex1Suff}
\hat{L}\geq \frac{1}{2}.
\end{equation}
When all $p^j(a_j)\ge 1/2$, then \eqref{eqn:Ex1Suff} and \eqref{eqn:4} automatically hold, but \eqref{eqn:Ex1Extra} is not necessarily satisfied. In that sense, the $K$-monotone, $K$-diminishing, and $K$-Go concavity conditions of Theorem~\ref{thm:myopicbounds3} 
are weaker sufficient conditions for achieving the bound $(1-(1-\frac{1}{K})^K)$  than the prefix monotone, diminishing-return, and postfix monotone conditions of Theorem~\ref{uniformstringmatroidthm1}.

\textbf{Adaptive Sensing}: Consider the situation where the additive noise set is independent but not identically distributed. Assume that $w_i$ is a Gaussian vector with mean zero and covariance $\sigma_i I$, where $I$ denotes the identity matrix.
Recall the problem formulation in Section~\ref{examples}. The objective function $f$ for this problem is as follows:
\[
f((B_1,\ldots, B_k))=\frac{1}{2}(\text{log det}(P_0)-\text{log det}(P_k)).
\]
where $P_0 = I$ and for $1\leq j\leq k-1$,
\[
P_j=\left(P_{j-1}^{-1}+\frac{1}{\sigma_j^2}B_j^TB_j\right)^{-1}.
\] 
From the expression above, it is easy to check that the order of $B_1,\ldots, B_k$ influences the objective function value under the assumption that $\sigma_1,\ldots, \sigma_k$ take different values. For example, \[f((A,B)) =  \frac{1}{2}\text{log det}\left(I+\frac{1}{\sigma_1^2}A^TA+\frac{1}{\sigma_2^2}B^TB\right)\] and \[f((B,A)) =  \frac{1}{2}\text{log det}\left(I+\frac{1}{\sigma_1^2}B^TB+\frac{1}{\sigma_2^2}A^TA\right).\] If $\sigma_1\neq \sigma_2$, then $f((A, B))\neq f((B,A))$.

By \cite{Liu2015}, some sufficient conditions for $f$ to be $K$-polymatroid and $K$-GO-concave are 
\begin{equation}
\label{suffcond1}
\sigma_{i+1}^2\geq \sigma_i^2
\end{equation}
for $i=1,\ldots, K-1$.

By \cite{Zhang2016}, to achieve the bound $(1-(1-1/K)^K)$, it requires both (\ref{suffcond1}) and 
\[\frac{b^{-2}}{a^{-2}-b^{-2}}\geq (K-1)^2(a^{-2}+b^{-2})+1,\]
where $[a,b]$ is the interval that contains all the $\sigma_i$'s. 

Comparing the sufficient conditions for achieving the same bound $(1-(1-1/K)^K)$ from \cite{Liu2015} and \cite{Zhang2016}, we see that the conditions from \cite{Liu2015} are weaker.

\section{Final Remarks}
\label{conclusions}
In this survey, we considered two classes of submodular maximization problems: set submodular maximization  and string submodular maximization. For set submodular optimization, we reviewed performance bounds for the greedy strategy under matroid constraints, improved performance bounds, and performance bounds for the batched greedy strategy. There are many important results about performance of the greedy strategy under some other constraints and conditions.  \cite{Wolsey1982}, \cite{SVIRIDENKO200441}, and  \cite{Kulik2009}  derived  performance bounds for the greedy strategy in submodular maximization problems subject to a knapsack constraint and multiple linear constraints. \cite{Bian2017} established performance bounds for the greedy strategy in monotone but nonsubmodular maximization problems under uniform matroid constraints.  People also investigated performance bounds for some variations of greedy strategies.  \cite{Calinescu2011} and \cite{Feldman2011} derived performance bounds for a randomized continuous greedy algorithm and a unified continuous greedy algorithm in monotone submodular maximization problems, respectively. \cite{Buchbinder2012}  established  performance  bounds for an adaptive greedy algorithm in unconstrained submodular maximization problems. They also derived 
 performance bounds for randomized greedy algorithms in nonmonotone submodular maximization problems  \citep{Buchbinder2014}.  \cite{mirzasoleiman16a} considered submodular maximization problems in a distributed fashion, and they derived  performance  bounds for  a two-stage greedy algorithm  under matroid or knapsack constraints. \cite{Qu2015} proposed a distributed greedy strategy and showed that it  has the same guarantee as the centralized greedy strategy.

 For string submodular optimization problems,  we reviewed performance bounds for the greedy strategy under matroid constraints. There are some related results on performance bounds for greedy strategies  in  string submodular maximization problems that were not reviewed in this paper. For example, \cite{Golovin2011} considered a particular class of partially observable adaptive stochastic optimization problems, and established performance bounds for the greedy strategy by introducing the notion of adaptive submodularity. \cite{Tschiatschek2017} derived performance bounds for a modified greedy strategy in submodular string optimization problems under uniform string matroid constraints.
%
%
%\cite{Buchbinder2012}

The scope of this study is limited to the performance of the greedy strategies in deterministic optimization problems where the objective function only involves actions. Potentially fruitful areas for further research include performance bounds for the greedy strategy in stochastic optimization problems, where the objective function involves states
and control actions, and real-world applications of the performance
bounds in the deterministic and stochastic settings.

\bibliographystyle{ametsoc2014}
\bibliography{bibfile} %note that this is a separate file

\end{document}